# Micro–macro population dynamics models of benthic algae with long-memory decay and generic growth

Hidekazu Yoshioka[1,*] and Kunihiko Hamagami[2]

[1] Graduate School of Advanced Science and Technology, Japan Advanced Institute of Science and Technology, 1-1 Asahidai, Nomi, Ishikawa, Japan (E-mail: yoshih@jaist.ac.jp)

[2] Faculty of Agriculture, Department of Food Production and Environmental Management, Iwate University, 3-18-8 Ueda, Morioka, Japan

[*] Corresponding author: yoshih@jaist.ac.jp, ORCID: 0000-0002-5293-3246

**Abstract**

Benthic algae as a primary producer in riverine ecosystems develop biofilms on the riverbed. Their population dynamics involve growth and decay processes, the former owing to the balance between biological proliferation and mortality, while the latter to mechanical abrasion because of the transport of sediment particles. Contrary to the assumptions of previous studies, the decay has experimentally been found to exhibit long-memory behavior, where the population decreases at an algebraic rate. However, the origin and mathematical theory of this phenomenon remain unresolved. The objective of this study is to introduce a novel mathematical model employing spin processes to describe microscopic biofilm dynamics. A spin process is a continuous-time jump process transitioning between states 0 and 1, and the continuum limit of these processes captures the long-memory decay and generates generic growth. The proposed framework leverages heterogeneous spin rates, achieved by appropriately superposing spin processes with distinct rates, to reproduce the long-memory decay. Computational simulations demonstrate the behavior of the model, particularly emphasizing rate-induced tipping phenomena. This mathematical model provides a computationally tractable interpretation of benthic algae dynamics and their long-term prediction, relevant to river-engineering applications.

***Statements & Declarations***

**Fundings** This study was supported by the Japan Society for the Promotion of Science (KAKENHI No. 25K07931).

**Competing interests** The authors have no relevant financial or non-financial interests to disclose.

**Acknowledgments** The authors gratefully acknowledge the helpful comments of three anonymous reviewers.

**Declaration of generative AI in scientific writing** The authors did not use generative AI for scientific writing of this manuscript.

**Data availability** The data will be made available upon reasonable request to the corresponding author.

**CRediT author statement (H. Yoshioka):** Conceptualization; Investigation; Methodology; Software; Validation; Visualization; Writing - original draft; Writing - review & editing. **(K. Hamagami):** Data curation; Funding acquisition; Project administration; Resources; Visualization; Writing - review & editing

## 1. Introduction

### 1.1 Research background

The sustainable coexistence of humans and the environment remains a critical concern, with water environments, such as rivers and lakes, being particularly impacted by anthropogenic activities. Examples include industrial water pollution [1], hydropeaking in dam-downstream regions causing fish stranding [2], lake area contraction due to irrigation water use [3], and the introduction and spread of invasive species [4].

Benthic algae, also known as attached algae, macroalgae, or periphyton, represent a key species in aquatic ecosystems as they function as primary producers driving food webs and nutrient cycling [5]. Consequently, the population dynamics of benthic algae are pivotal in evaluating the sustainability of human-environment interactions. Eutrophication owing to excess nutrient input often results in the overgrowth of benthic algae as stable equilibria [6,7,8], leading to secondary destabilizing effects on food webs such as the dispersal of invasive snails [9] and disruptions to the aquatic carbon balance [10] (see also **Photo 1**). Monitoring and regulating benthic algae populations are thus essential components of aquatic environmental management for sustainability.

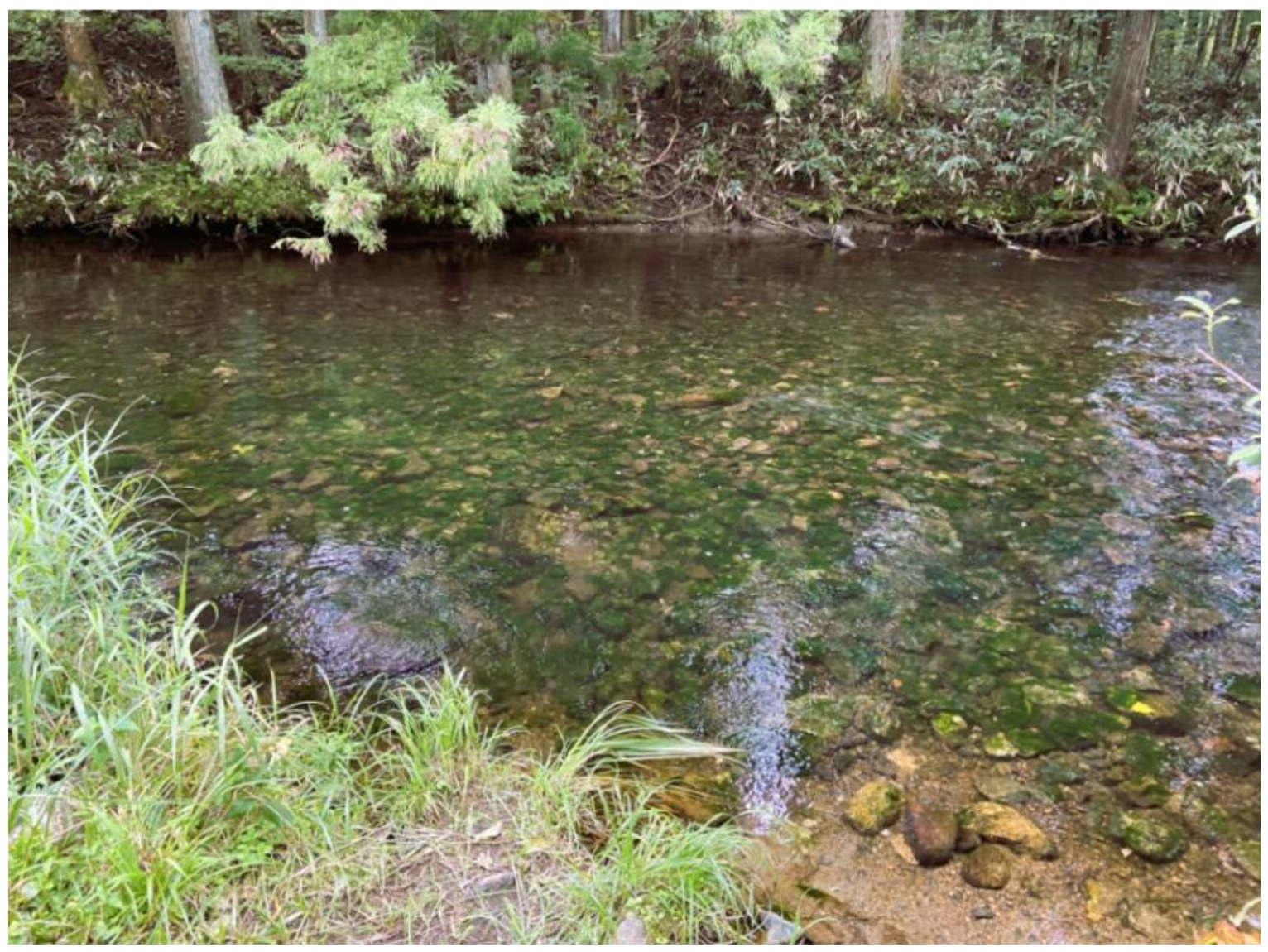

**Photo 1. Benthic** algae bloom at Nakatsu River in Iwate Prefecture, Japan. This photo was taken by the second author (Kunihiko Hamagami) on September 12, 2021.

The population dynamics of benthic algae encompass both growth and decay processes. Growth arises from development regulated by chemical and biological factors and is typically modeled using ordinary differential equations (ODEs) or discrete-time versions, such as logistic models [11] and Allee models [12]. Coefficients in these models are primarily influenced by nutrient availability [13,14], water flow velocity and turbulence [15,16,17], photosynthesis [18,19], and their interactions [20,21]. Lima et al. [22] investigated hydrological variables impacting river biofilms. The mathematical aspects of growth dynamics have been elaborated by incorporating species interactions and spatial distribution [23,24,25] and persistence [26]. As reviewed above, ODEs and related models have extensively been used as fundamental tools for studying population dynamics.

Conversely, the decay in algae population dynamics is attributed to abrasion by the transport of sediment (sand and gravel) particles. Specifically, the removal of benthic algae on riverbeds results from collision with sediment particles, recognized as the primary mechanism governing algae population dynamics [27]. In natural river environments, abrasion is induced when flow discharge or velocity exceeds some threshold [28,29], suggesting that frequent flood pulses effectively regulate benthic algae populations.

The abrasion effects have conventionally been incorporated into population dynamics models of benthic algae by assuming exponential decay, where the decay rate depends on flow conditions as reported in applied and modeling studies [30,31,32,33]; however, recent experimental findings have suggested that the population decay of benthic algae due to abrasion is not exponential ( $\exp(-t)$ , where $t$ denotes time with appropriate scaling) but can be algebraic ( $t^{-\alpha}$ with $\alpha > 0$ ) [34,35]. **Figure 1** illustrates experimental data showing population dynamics measured via the algae surface coverage ratio of a hemisphere of fixed radius in sediment-laden water flow. The data align closely with an algebraic, long-memory decay curve, while the exponential model fails to capture the dynamics, underestimating the decay at a short timescale and overestimating it at a long timescale. This implies that a model with nonexponential decay would be necessary for advanced modeling of the population dynamics of benthic algae.

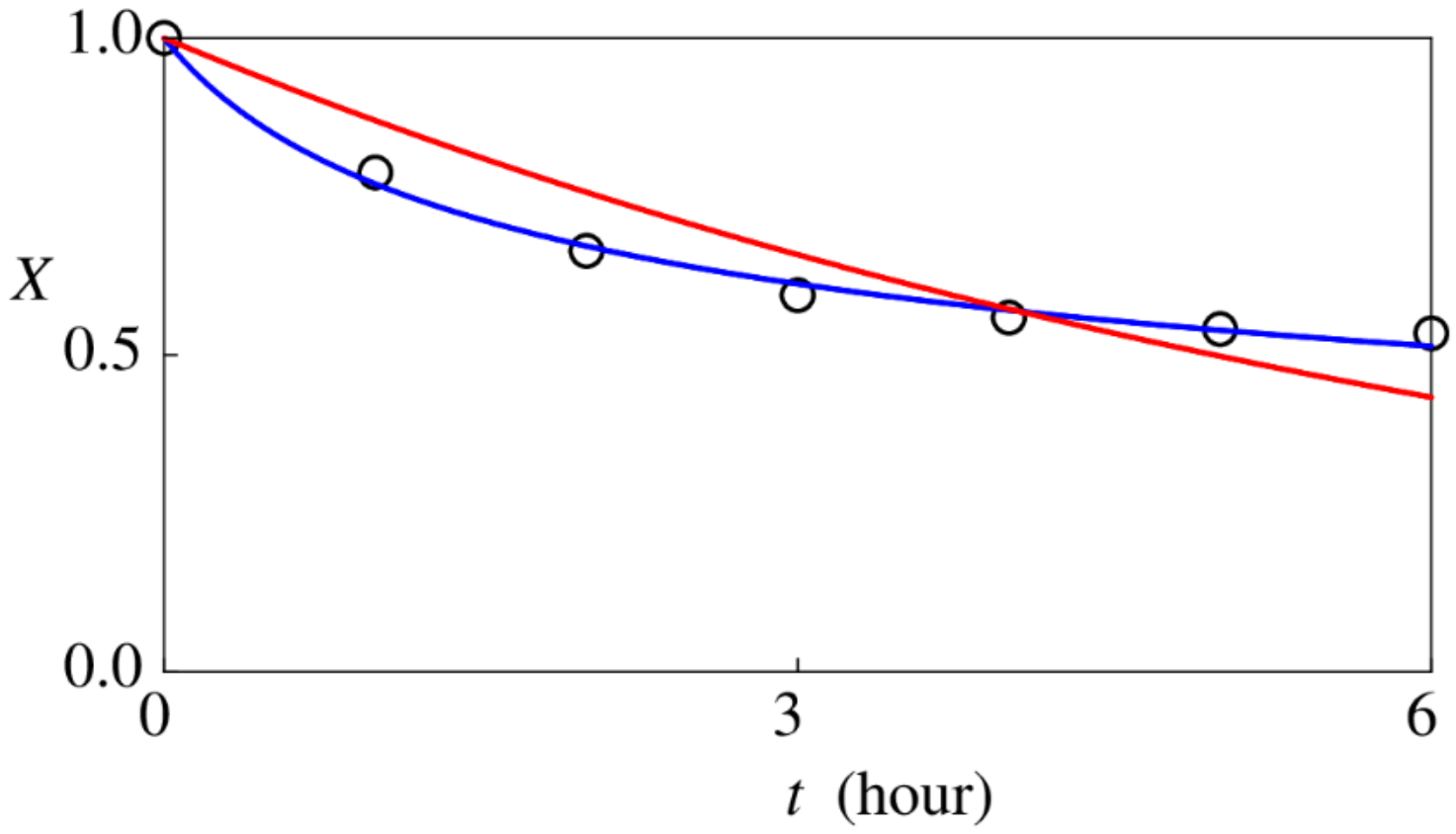

**Figure 1.** Decay of benthic algae population represented by the algae covering ratio on the surface of a hemisphere with a specified radius in sediment-laden water flow. Black circles denote data points. The blue curve depicts the long-memory fit, while the red curve illustrates the exponential fit. The experimental conditions are detailed in **Section 3**. The long-memory fit is better, as shown in the figure.

Theoretically, the difference between the exponential and algebraic decays is important because it implies that human interventions to restore degraded riverbed environments, such as sediment replenishment [36,37], may not work effectively because decay may not dominate growth (Proposition 3 in Yoshioka and Hamagami [35]). Mathematical models describing the population dynamics of benthic algae, which incorporate both growth and abrasion, should therefore account for the algebraic decay due to abrasion. However, the theoretical framework for this phenomenon remains poorly understood, motivating this study. The key question of this study is therefore as follows: How can the population dynamics with the long-memory decay of benthic algae be modeled? Theoretical descriptions of such a phenomenon would serve as a potential new step toward modeling the population dynamics of benthic algae reviewed above.

### 1.2 Aim and contribution

The objective of this study is to develop a novel theoretical framework capable of representing the nonexponential population decay of benthic algae. Furthermore, it seeks to establish a comprehensive theory that integrates growth and decay dynamics, enabling applications to both long- and short-term problems. The proposed model offers a novel interpretation of the dynamics of benthic algae, providing a computationally feasible interpretation. See also **Figure 2**.

**Previous studies**

- ✓ Benthic algae population decays exponentially.
- ✓ Tractable due to the exponential nature.
- ✓ The exponential nature is not consistent with the experiment.

**This study**

- ✓ Benthic algae population does not necessarily decay exponentially.
- ✓ Mathematical foundation of the model is developed.
- ✓ The model is experimentally validated

**Figure 2.** Comparison diagram between previous studies and this study.

Our contributions are elaborated below. Our theory derives macroscopic population dynamics from microscopic dynamics, assuming that, for a given nonflat riverbed area, there exists a variety of decaying timescales, i.e., heterogeneity in the exponential decay rate. The reason that heterogeneous exponential decay leads to long-memory decay can be conceptually understood through the following formula [34]:

$$\underbrace{\int_0^{+\infty}\underbrace{\left(\frac{1}{\Gamma(\alpha)}\frac{R^{\alpha-1}}{\beta^{\alpha}}e^{-\frac{R}{\beta}}\right)}_{\text{Heterogeneity}}\underbrace{e^{-Rt}}_{\text{Exponential decay}}\mathrm{d}R}_{\text{Superposition}} = \underbrace{\frac{1}{(1+\beta t)^{\alpha}}}_{\text{Long-memory decay}}, \quad t \geq 0, \tag{1}$$

where $\alpha, \beta > 0$ are constants and $\Gamma$ is the classical gamma function. Considering (1), the emergence of long-memory decay (the right-hand side) is attributed to the superposition of the heterogeneous exponential decay rates (the left-hand side) under a gamma-type probability measure: $\sim R^{\alpha-1}e^{-R/\beta}\mathrm{d}R$. This study investigates the mechanism underlying (1) and its generalization. The key advantage of the present framework based on superposition is its theoretical simplicity; it does not need a modeling framework that is not local in time and often needs a sophisticated notion such as fractional motions [38]. In contrast, our superposition approach is just a summation of independent processes in the simplest case and is therefore theoretically simpler than these existing models.

This paper uses the spin process (Aalen [39], Chapter 7 in Capasso and Bakstein [40]), alternatively referred to as a hazard process, which is a continuous-time stochastic process that transitions randomly between states 1 and 0. Here, state 1 represents the presence of benthic algae within an infinitesimally small segment of the riverbed, whereas state 0 indicates its absence. The spin processes examined in this work are heterogeneous, characterized by distinct spin rates across individual processes. This paper demonstrates that the nonexponential decay arises by superposing the heterogeneous spin processes because of weak convergence (Chapter 13 in Klenke [41]). This contribution substantiates the conjecture proposed by Yoshioka et al. [34], which used a more complex, doubly stochastic framework. Previous bottom-up approaches have addressed simulations of macroscopic population dynamics based on microscopic descriptions, such as multifractal analysis of observed periphyton biomass on a plate [42] and multispecies probabilistic individual-based models [43], but they did not explicitly consider population dynamics accounting for both growth and decay.

Note that spin processes and their reversible versions have extensively been studied across various domains, including survival analysis [44], clinical trials [45], health and life insurance [46], option pricing and optimal contracts [47,48], machine maintenance [49], and moral hazard problems [50]. Mathematical investigations of spin processes have also been conducted in general settings [51,52]. Hence, this paper explores a new use of spin processes, modeling the population dynamics of benthic algae and extending the application area of these processes.

A key advantage of the proposed theory is its ability to incorporate a broad class of growth functions, including classical logistic and Allee-type dynamics. The framework enables a unified treatment of both homogeneous and heterogeneous growth dynamics of benthic algae by modeling reversible spin processes. Furthermore, the flexibility of the theory permits the inclusion of time-dependent parameters in the model coefficients, facilitating the investigation of population dynamics under environmental fluctuations that may lead to rate-induced tipping [53,54]. The theoretical simplicity of the proposed model in reproducing macroscopic population dynamics from microscopic stochastic processes potentially applies to modeling other species dynamics, including macroalga and both aquatic and nonaquatic vegetation, which would be seen worldwide [55,56,57,58,59].

The governing equation of the population dynamics model in this study is not an ODE but rather a continuum of ODEs, formulated as an integro-differential equation (IDE). The unique existence of its solution is established in the space of bounded and integrable functions under appropriate assumptions on the model coefficients, with the solution satisfying the model in a classical pointwise sense. The inherent nonlocality, which induces long-memory decay, distinguishes the proposed model from conventional ODE-based population dynamics models. Furthermore, the model construction based on spin processes provides a framework for simulating population dynamics and quantifying model uncertainty arising from the misspecification of heterogeneity. This paper computationally applies the population dynamics model to a rate-tipping problem and an advanced scenario. Tipping phenomena are explored for microscopic, intermediate, and macroscopic models in a consistent manner. This study thus contributes to the formulation, mathematical analysis, and computational application of a novel population dynamics model for benthic algae.

The theoretical results obtained indicate that the mathematical model presented in this paper is mathematically well-posed (having a unique global solution: **Proposition 4**) and stable against disturbances in the input data (**Proposition 5**). In other words, it is a model that can withstand engineering applications where precise observational data is not always available. The experimental results (e.g., **Figure 1**) suggest that assuming an exponential rate of decay in attached algae leads to an underestimation of the duration of algal growth. The mathematical model presented in this paper has the advantage of being able to compensate for decay that is far slower than exponential. This is likely to become important when considering environmental management policy design and costs in practical civil engineering projects.

The remainder of this study is organized as follows. **Section 2** formulates the spin processes governing microscopic biofilm dynamics. The population dynamics model, the main outcome of this study, is developed in this section. **Section 3** focuses on computational analysis based on the population dynamics model. **Section 4** provides a summary of the findings and outlines future directions for this research. **Appendix** contains proofs and supplementary results.

## 2. Mathematical model

The streamline of our mathematical modeling is as follows. First, spin processes that play a role in this paper are introduced along with their macroscopic, continuum limit (**Sections 2.1-2.2**). Second, the growth dynamics of benthic algae are studied in view of spin processes (**Section 2.3**). Third, the population dynamics model of benthic algae is introduced based on spin processes and its well-posedness is investigated (**Section 2.4**). **Table 1** summarizes key notations.

**Table 1.** Key notations.

| Notation | Explanation |
|---|---|
| $t$ | Time |
| $D$ | Physical domain |
| $D_i$ | $i$ th subdomain of $D$ |
| $M \in \mathbb{N}$ | Total number of subdomains $D_i$ |
| $R > 0$ | Distributed spin rate (decay rate) |
| $G$ | Growth coefficient |
| $F$ | Probability distribution of $R$ |
| $x_{i,t} \in \mathbb{N}$ | Spin on $D_i$ at time $t$ in a discrete $R$ setting |
| $X_t^{(M)}$ | Total population at time $t$ in a discrete $R$ setting |
| $x_t(R)$ | Population density with the decay rate $R$ at time $t$ in a continuum limit |
| $X_t$ | Total population at time $t$ in a continuum limit |

This paper considers the population dynamics of benthic algae within a defined riverbed domain $D$, with an area denoted by $|D| > 0$. The carrying capacity of the benthic algae in $D$ is normalized to 1, such that the population in $D$ is a time-dependent variable constrained between the minimum 0 and maximum 1. In applications, the population can be quantified either as biomass or a coverage ratio of $D$. Using the latter approach, the population equals 1 if $D$ is completely covered by benthic algae and 0 if there are no algae in $D$. The latter approach is used here because it aligns with our experimental framework described subsequently. For simplicity, $D$ is assumed to be an isolated habitat, precluding any exchange of population with external environments, as our primary focus is on decay and growth.

This paper is based on a complete probability space $(\Omega, \mathbb{F}, \mathbb{P})$, where $\Omega$ represents the collection of all events, $\mathbb{F}$ is a filtration, and $\mathbb{P}$ is a probability function as conventionally defined in stochastic models (Definition 1.1 in Capasso and Bakstein [40]). Expectation is denoted by $\mathbb{E}$, and variance is denoted by $\mathbb{V}$. Time $t$ is a nonnegative parameter. The set of càdlàg (right-continuous with left limit) processes over the time interval $[0,T]$ with the range $\Xi$ is denoted as $\mathbb{D}([0,T],\Xi)$. The left limit of a stochastic process $x$ at time $t$ is expressed as $x_{t-}$.

### 2.1 Spin process

A fundamental spin process is explained. Let $x = (x_t)_{t \geq 0}$ represent an irreversible spin process that spins from 1 to 0, modeled as a counting process with the initial condition $x_0 = 1$ such that $x_t = \mathbb{I}(t < \tau)$, where $\mathbb{I}(a < b)$ is the indicator function such that $\mathbb{I}(a < b) = 1$ if $a < b$ and $\mathbb{I}(a < b) = 0$ otherwise, and $\tau$ is a random variable that follows the exponential distribution with the mean $R^{-1}$ where $R > 0$ is the spin rate. The spin process is càdlàg. Furthermore, the following process constitutes an orthogonal and square-integrable martingale concerning a natural filtration generated by $x$ (Lemma 4.1 in Aalen [39]):

$$x_t - \int_0^t R(1 - x_{s-}) \mathrm{d}s, \quad t \geq 0. \tag{2}$$

An elementary calculation illustrates that the average of $x_t$ is $e^{-Rt}$ and the variance is

$$\mathbb{V}[x_t] = e^{-Rt}\left(1 - e^{-Rt}\right), \quad t \geq 0. \tag{3}$$

Note that $x_t^2 = x_t$.

### 2.2 Population decay

This subsection derives the long-memory decay of the algae population under the assumption that there exists a continuous flow of sediment particles that potentially removes the algae from the riverbed.

This paper adopts a model where a continuum of spin processes, characterized by distinct spin rates, exists within the domain $D$. The heterogeneity in natural settings arises from the nonflat riverbed, which induces spatially nonuniform local flow conditions in $D$. In laterally homogeneous experimental settings with a channel bed featuring semicircular bumps, Hamagami et al. [60] demonstrated that the local water flow along the longitudinal direction in the channel exhibits heterogeneity. This heterogeneity results in varying shear stresses on the upstream and downstream surfaces of bumps, identified as $D$ in our framework, which in turn causes spatially heterogeneous decay of benthic algae. Hamagami et al. [60] found that the detachment process of benthic algae is spatially heterogeneous, with more pronounced detachment near the top of the bumps. **Figure 1** illustrates the results from additional laboratory experiments on the decay of benthic algae on hemispheres, where heterogeneous water flow conditions have been observed [61]. In this scenario, detachment is more significant on surfaces oriented toward the flow current, which is due to the spatial flow difference between the upstream and downstream sides of the hemispheres. Homogeneous models corresponding to exponential decay do not apply in this case, motivating the use of a heterogeneous model.

The stochastic system that serves as a microscopic representation of algae population dynamics within the domain $D$ is introduced. The domain $D$ is divided into $M \in \mathbb{N}$ subdomains $\{D_i\}_{i=1,2,3,\ldots,M}$. Each $D_i$ is associated with an irreversible spin process that transitions from 1 to 0 at a spin rate $R_i > 0$. Without loss of generality, $D_i$ is assumed to have an equal area of $M^{-1}|D|$. Additionally, the pointwise spin rate $R$ in $D$ is assumed to be distributed according to a probability measure $F = F(\mathrm{d}R)$. This implies

$$F\left(\underline{R},\overline{R}\right)=\frac{\left|D\left(\underline{R},\overline{R}\right)\right|}{\left|D\right|}. \tag{4}$$

Here, $D\left(\underline{R},\overline{R}\right)$ denotes the union of subdomains of $D$ such that the spin rate lies within the range $\left(\underline{R},\overline{R}\right)$ with $\underline{R}<\overline{R}$; $D\left(0,+\infty\right)=D$. Here, the discretization of $D$ adopts the quantile-based method utilized in prior investigations of Markovian lifts (e.g., Yoshioka [62]). It defines the quantile $R_i$:

$$\int_0^{R_i} F\left(\mathrm{d}r\right)=\frac{2i-1}{2M},\quad i=1,2,3,...,M. \tag{5}$$

Consequently, the discretized version $F_M$ of $F$ is set as

$$F_M\left(\mathrm{d}r\right)=\frac{1}{M}\sum_{i=1}^{M}\delta\left(R-R_i\right),\quad R>0. \tag{6}$$

This discretization is useful because of a uniform bound between $F$ and $F_M$ (see **Proof of Proposition 1**).

The population in $D_i$ is assumed to follow an irreversible spin process $x_i=\left(x_{i,t}\right)_{t\geq 0}$ that spins from 1 to 0, with each $x_i$ being mutually independent. This study operates under the following technical assumption that $F$ is nonsingular, which is not critical in our application where a gamma-type distribution is assumed.

***Assumption 1*** *The probability measure $F$ admits a density.*

Physically, this assumption represents the situation where the decay rate is gradually distributed over the domain.

For each $M$, the algae population $X^{(M)}=\left(X_t^{(M)}\right)_{t\geq 0}$ in $D$ is expressed as follows:

$$X_t^{(M)}=\frac{1}{M}\sum_{i=1}^{M}x_{i,t},\quad t\geq 0. \tag{7}$$

The following **Proposition 1** is the first primary result in this study and generalizes (1).

***Proposition 1***

*It follows in the sense of probability that*

$$\lim_{M\to+\infty} X_t^{(M)}=\int_0^{+\infty}e^{-Rt}F\left(\mathrm{d}R\right),\quad t\geq 0. \tag{8}$$

The proof is based on the derivation of the average and variance of spin processes, and how the macroscopic limit will arise is computationally demonstrated in **Section 4**.

***Remark 1*** One may also consider a continuum of spin processes by employing the Fubini extension, which constitutes a proper mathematical framework for a continuum of independent stochastic processes [63]. However, this approach involves significant complexity, requiring nontrivial integration on two measure spaces [64].

***Remark 2*** This study primarily examines the limit $M \to +\infty$; however, finite $M$ cases may also warrant investigation if benthic algae population dynamics span patchy habitats characterized by a piecewise constant environment.

### 2.3 Population growth

This subsection considers the growth dynamics associated with the following growth rate generalizing classical logistic and Allee ones:

$$G(x) = rx(1-x)g(x), \quad x \in \mathbb{R}, \tag{9}$$

where $r > 0$ is the intrinsic growth rate that is assumed to be 1 unless otherwise specified and without loss of any generality, $g : \mathbb{R} \to \mathbb{R}$ is assumed to be decomposed into the positive and negative parts:

$$g(x) = g^{+}(x) - g^{-}(x), \quad x \in \mathbb{R}. \tag{10}$$

Here, $g^{+}, g^{-} : \mathbb{R} \to [0, +\infty)$ represent nonnegative, bounded, and Lipschitz continuous functions (see **Remark 3**). Then, the following ODE admits a unique nonnegative solution (Proposition 4.2 in Magnus [65]) that is bounded between 0 and 1:

$$\frac{\mathrm{d}X_t}{\mathrm{d}t} = G(X_t), \quad t > 0 \tag{11}$$

subject to an initial condition $X_0 \in [0,1]$. The positive component $g^{+}$ denotes standard growth functions such as logistic growth, while the negative component $g^{-}$ reflects the Allee effect, where the growth rate becomes negative; the latter case chooses $g^{+}(x) = x$ and $g^{-}(x) = a$ (with some regularization) with a constant $a \in (0,1)$. The logistic case follows by setting $g^{+}(x) = 1$ and $g^{-}(x) = 0$.

***Remark 3*** The assumption that $g^{+}, g^{-}$ are nonnegative and bounded imposes excessive for the Allee-type one. The regularized version, such as $g^{+}(x) = \max\{0, \min\{1, x\}\}$, may be restrictive; however, the solution of the proposed population dynamics model is inherently bounded between 0 and 1, rendering regularization such as the aforementioned $g^{+}$ inactive (see **Proof of Proposition 2**).

It is assumed that each spin process, again denoted by $\{x_{i,t}\}_{i=1,2,3,\ldots,M}$ ($t \geq 0$), is now reversible:

$$x_{i,t} = x_{i,0} + \int_0^t \mathbb{I}(x_{i,s-} = 0)\,\mathrm{d}N_{i,s} - \int_0^t \mathbb{I}(x_{i,s-} = 1)\,\mathrm{d}L_{i,s}, \quad i = 1,2,3,\ldots,M, \tag{12}$$

where $\{N_{i,t}\}_{i=1,2,3,...,M}$ represents a series of point processes with a common jump intensity of $X_{t-}^{(M)} g^{+}\left(X_{t-}^{(M)}\right)$, and $\{L_{i,t}\}_{i=1,2,3,...,M}$ represents another series of point processes with a common jump intensity of $\left(1-X_{t-}^{(M)}\right) g^{-}\left(X_{t-}^{(M)}\right)$, where the notation (7) is used again. More rigorously, $N_{i,t}$ (resp., $L_{i,t}$) is expressed via some Poisson random measure $\bar{N}_i(\mathrm{d}u\mathrm{d}z\mathrm{d}t)$ (resp., $\bar{L}_i(\mathrm{d}u\mathrm{d}z\mathrm{d}t)$) on $(0,+\infty)^3$:

$$\mathrm{d}N_{i,t} = \int_{u=0}^{u=X_{t-}^{(M)} g^{+}\left(X_{t-}^{(M)}\right)} \int_{z=0}^{z=1} \bar{N}_i(\mathrm{d}u\mathrm{d}z\mathrm{d}t) \text{ (resp., } \mathrm{d}L_{i,t} = \int_{u=0}^{u=\left(1-X_{t-}^{(M)}\right) g^{-}\left(X_{t-}^{(M)}\right)} \int_{z=0}^{z=1} \bar{L}_i(\mathrm{d}u\mathrm{d}z\mathrm{d}t) \text{ ).} \quad (13)$$

It is also assumed that $\bar{N}_i$ and $\bar{L}_i$ are mutually independent. The initial condition $x_{i,0}$ is either 0 or 1, and hence, $x_{i,t}$ is also because jumps of $N_i$ and $L_i$ are 1. The stochastic system (12) admits at most one càdlàg pathwise solution that is bounded between 0 and 1 (see **Lemma 1** in **Appendix**).

Because $\mathbb{I}(x_{i,s-}=1)=x_{i,s-}$ and $\mathbb{I}(x_{i,s-}=0)=1-x_{i,s-}$, (12) is rewritten as follows so that the martingale component (mean-zero fluctuation part) is explicitly found:

$$\begin{aligned} & x_{i,t}-x_{i,0}-\int_0^t\left(1-x_{i,s}\right) X_s^{(M)} g^{+}\left(X_s^{(M)}\right) \mathrm{d}s+\int_0^t x_{i,s}\left(1-X_s^{(M)}\right) g^{-}\left(X_s^{(M)}\right) \mathrm{d}s \\ & =\int_0^t\left(1-x_{i,s-}\right) \mathrm{d}\tilde{N}_{i,s}-\int_0^t x_{i,s-} \mathrm{d}\tilde{L}_{i,s} \end{aligned} \quad (14)$$

with compensated processes $\mathrm{d}\tilde{N}_{i,s}=\mathrm{d}N_{i,s}-X_s^{(M)} g^{+}\left(X_s^{(M)}\right)\mathrm{d}s$ and $\mathrm{d}\tilde{L}_{i,s}=\mathrm{d}L_{i,s}-\left(1-X_s^{(M)}\right) g^{-}\left(X_s^{(M)}\right)\mathrm{d}s$. The right-hand side of (14) is the martingale part, while the left-hand side is its alternative representation and is a key for deriving the growth model. Summing up (14) for $i=1,2,3,...,M$ along with (10) yields

$$X_t^{(M)}-X_0^{(M)}-\int_0^t\left(1-X_s^{(M)}\right) X_s^{(M)} g\left(X_s^{(M)}\right) \mathrm{d}s=\frac{1}{M} \sum_{i=1}^M\left\{\int_0^t\left(1-x_{i,s-}\right) \mathrm{d}\tilde{N}_{i,s}-\int_0^t x_{i,s-} \mathrm{d}\tilde{L}_{i,s}\right\}. \quad (15)$$

The following **Proposition 2** shows that the right-hand side of (15) vanishes as $M \rightarrow +\infty$.

***Proposition 2***

*Fix* $T>0$. *For each* $t \in [0,T]$, *it follows that*

$$\lim_{t \rightarrow +\infty} \mathbb{E}\left[\left(\frac{1}{M} \sum_{i=1}^M\left\{\int_0^t\left(1-x_{i,s-}\right) \mathrm{d}\tilde{N}_{i,s}-\int_0^t x_{i,s-} \mathrm{d}\tilde{L}_{i,s}\right\}\right)^2\right]=0 \quad (16)$$

*and hence*

$$\lim_{t \rightarrow +\infty} \mathbb{E}\left[\left(X_t^{(M)}-X_0^{(M)}-\int_0^t\left(1-X_s^{(M)}\right) X_s^{(M)} g\left(X_s^{(M)}\right) \mathrm{d}s\right)^2\right]=0 . \quad (17)$$

According to **Proposition 2**, the limit $\hat{X}_t=\lim_{t \rightarrow +\infty} X_t^{(M)}$ should satisfy

$$\hat{X}_t-\hat{X}_0-\int_0^t\left(1-\hat{X}_s\right) \hat{X}_s g\left(\hat{X}_s\right) \mathrm{d}s=0, \quad (18)$$

which is an integrated version of (11). This statement is justified in the following **Proposition 3** (Chapter 3.1 in Bansaye and Méléard [66]).

***Proposition 3***

*Assume that* $\lim_{M\to+\infty} X_0^{(M)} = \hat{X}_0$ *in the sense of probability where* $\hat{X}_0 \in [0,1]$ *is a constant. Fix* $T > 0$. *For each* $t \in [0,T]$, *the equation (18) holds true in the sense of probability.*

***Remark 4*** Our approach also applies to a more generic growth model such as $G(x) = (1-x)g(x)$ (there is a generation of the population when $G(0) = g(0) > 0$) if the corresponding solution $\hat{X}$ to the ODE (11), subject to an initial condition $\hat{X}_0 \in [0,1]$, is still bounded in $[0,1]$ globally in time.

## 2.4 Population dynamics model

### 2.4.1 Formulation

Now, a population dynamics model that integrates both growth and decay dynamics is examined. It is assumed that each spin process, again denoted by $\{x_{i,t}\}_{i=1,2,3,\ldots,M}$ ($t \geq 0$), satisfies

$$x_{i,t} = x_{i,0} + \int_0^t \left(1 - x_{i,s-}\right) \mathrm{d}N_{i,s} - \int_0^t x_{i,s-} \mathrm{d}L_{i,s} - \int_0^t x_{i,s-} \mathrm{d}K_{i,s}, \quad i = 1,2,3,\ldots,M, \tag{19}$$

where each of $N_i$ and $L_i$ are defined as in the previous subsection and the sequence $\{K_{i,t}\}_{i=1,2,3,\ldots,M}$ is a series of point processes where $K_{i,t}$ has the jump rate $R_i$. They are defined by some Poisson random measures, as in (13), that are mutually independent. Therefore, the right-hand side of (19) represents the initial condition, state-dependent positive growth, state-dependent negative growth, and decay. The model assumes that the growth and decay dynamics are modeled by different jump processes, assuming that they are driven by physically different events (abrasion and proliferation).

Under the limit $M \to +\infty$, the process $X^{(M)}$ in this case should converge to some deterministic process $\hat{X}$ such that

$$\hat{X}_t = \int_0^{+\infty} \hat{x}_t(R) F(\mathrm{d}R), \quad t \geq 0 \tag{20}$$

along with

$$\frac{\mathrm{d}\hat{x}_t(R)}{\mathrm{d}t} = \hat{X}_t g^+\left(\hat{X}_t\right) - \left\{R + \hat{X}_t g^+\left(\hat{X}_t\right) + \left(1 - \hat{X}_t\right) g^-\left(\hat{X}_t\right)\right\} \hat{x}_t(R), \quad R > 0 \tag{21}$$

subject to an initial condition $\hat{x}_0(R) \in [0,1]$ ($R > 0$). The time evolution of $\hat{X}$ is given by

$$\frac{\mathrm{d}\hat{X}_t(R)}{\mathrm{d}t} = \underbrace{\hat{X}_t g\left(\hat{X}_t\right)}_{\text{Growth}} - \underbrace{\int_0^{+\infty} R\hat{x}_t(R) F(\mathrm{d}R)}_{\text{Decay}}, \quad t \geq 0. \tag{22}$$

In this view, the population $\hat{X}$ solves (20). As shown on the right-hand side of (22), the long-memory decay of the population appears as an integral term that accounts for the distributed decay rate.

### 2.4.2 Well-posedness of integro-differential equation

It is shown that the IDE (21) admits a mild solution that is bounded and continuous and that the solution is actually a classical solution, i.e., it satisfies the equation pointwise.

The Banach space $\mathbb{L}_1$ is set as

$$\mathbb{L}_1 = \left\{ \hat{x} : (0,+\infty) \to \mathbb{R};\ \|\hat{x}\| = \int_0^{+\infty} |\hat{x}(R)| F(\mathrm{d}R) < +\infty \right\}. \tag{23}$$

Another Banach space based on $\mathbb{L}_1$ is also introduced: for each $T > 0$,

$$\mathbb{L}_{1,T} = \left\{ \hat{x} : [0,T] \times (0,+\infty) \to \mathbb{R};\ \sup_{0 \le t \le T} \|\hat{x}_t\| < +\infty \right\}. \tag{24}$$

The target of this section is the following nonlinear equation:

$$\begin{aligned} \hat{x}_t(R) &= e^{-Rt}\hat{x}_0(R) + \int_0^t e^{-R(t-s)} \left\{ \tilde{X}_s g^+(\tilde{X}_s) - \left\{ \tilde{X}_s g^+(\tilde{X}_s) + (1-\tilde{X}_s) g^-(\tilde{X}_s) \right\} \tilde{x}_s(R) \right\} \mathrm{d}s \\ &\left(= \mathbb{G}[\hat{x}](t,R)\right) \end{aligned},\quad t \ge 0,\quad R > 0, \tag{25}$$

where

$$\tilde{x}_s(R) = \max\left\{0, \min\left\{1, \hat{x}_s(R)\right\}\right\} \text{ and } \tilde{X}_s = \int_0^{+\infty} \tilde{x}_s(R) F(\mathrm{d}R), \tag{26}$$

and $\mathbb{G}$ for a measurable function $\hat{x}$ is a mapping from $[0,+\infty) \times (0,+\infty)$ to $\mathbb{R}$. The nonlinear equation (25) is a mild version of the regularized model

$$\frac{\mathrm{d}\hat{x}_t(R)}{\mathrm{d}t} = \tilde{X}_t g^+(\tilde{X}_t) - \left\{ R + \tilde{X}_t g^+(\tilde{X}_t) + (1-\tilde{X}_t) g^-(\tilde{X}_t) \right\} \tilde{x}_t(R),\quad R > 0 \tag{27}$$

subject to $\hat{x}_0(R) \in [0,1]$ ($R > 0$). The following proposition establishes the well-posedness and regularity of the IDE (21).

***Proposition 4***

*Assume that $\|\hat{x}_0\| < +\infty$. The regularized model (27) admits a unique solution in $\mathbb{L}_{1,T}$ for any $T > 0$. The solution is continuously differentiable at each $t > 0$, and hence it satisfies the equation (27) point-wise. The range of the solution is at most $[0,1]$, and hence the solution also satisfies the IDE (21). Moreover, the IDE (21) admits a unique solution in $\mathbb{L}_{1,T}$ for any $T > 0$.*

In applications, identifying coefficients, parameter values, and initial conditions of the population dynamics model without any error is often infeasible; thus, analyzing the impact of modeling errors on the solution constitutes a fundamental task. The proposition below establishes the influence of errors in the probability measure $F$ and initial conditions of the population dynamics model.

***Proposition 5***

*Let* $\hat{x}, \hat{y}$ *be the solutions to the IDE (21) with the initial conditions* $\hat{x}_0, \hat{y}_0$ *and probability measures* $F = F_x, F_y$*, respectively. Then, there exists a constant* $C > 0$ *independent from* $\hat{x}, \hat{y}$ *such that*

$$\sup_{R>0} \left| \hat{x}_t(R) - \hat{y}_t(R) \right| \le \left( \sup_{R>0} \left| \hat{x}_0(R) - \hat{y}_0(R) \right| + Ct \left\| F_x - F_y \right\|_{\mathrm{TV}} \right) \int_0^t e^{Cs} \mathrm{d}s, \quad t > 0. \tag{28}$$

*Here,* $\left\| F_x - F_y \right\|_{\mathrm{TV}}$ *is the total variation norm between* $F_x, F_y$ *given by*

$$\left\| F_x - F_y \right\|_{\mathrm{TV}} = \frac{1}{2} \sup_{\varphi} \left| \int_0^{+\infty} \varphi(R) \left( F_x(\mathrm{d}R) - F_y(\mathrm{d}R) \right) \right|, \tag{29}$$

*and the supremum in (29) is taken with respect to measurable functions* $\varphi : (0, +\infty) \to \mathbb{R}$ *bounded between 0 and 1.*

***Remark 5*** Population dynamics models employed in various disciplines, including cellular and bacterial dynamics [67,68], evolutionary games [69], and voting dynamics [70,71], originate from microscopic stochastic processes. These models predominantly assume homogeneous agents (where $F$ is represented by a Dirac delta in our framework), while ours assumes their heterogeneity. The proposed approach thus extends these models to incorporate heterogeneous agent-based systems [72].

## 3. Computational analysis

### 3.1 Experimental setting

Two hydraulic experiments were conducted using a flume located at the Faculty of Agriculture, Iwate University, Japan (**Figure 3**). While the experimental setup has been detailed in previous studies [34,35], it is reiterated here for completeness of this study along with a relevant open data set (**Section A.1** in **Appendix**).

In each experiment, a steady flow with a prescribed flow discharge was established after the placement of hemispherical structures mimicking boulders on the flume bed. These hemispheres have been coated with the green benthic alga *Cladophora glomerata Kützing*. Once positioned, sand particles were continuously added at the upstream end of the flume to create a steady sediment-laden water flow. Over time, the algae cover decayed, and the covering ratio (the percentage of the surface area of the hemisphere covered by the algae) of each hemisphere was measured every hour for 6 hours. To quantify the detachment of benthic algae in relation to sediment transport, a visual evaluation method was employed. A grid was overlaid on each hemispherical surface, and photographs were taken from four directions: front, back, left, and right. For each direction, the proportion of grid cells covered by algae was determined. These proportions were then averaged across the four directions to provide an overall index of detachment extent. This method was adapted from previous laboratory studies using artificial riverbeds with hemispherical objects [61]. The two experimental conditions are summarized in **Table 2**. These experiments were conducted under supercritical flow conditions, characterized by a Froude number between 1.65 and 1.69, indicative of steep stream flows with continuous sediment transport, facilitating continuous removal of benthic algae from the riverbed.

The domain $D$ is defined as the surface of the union of hemispheres depicted in **Figure 3**, corresponding to the experimental setup. For each experiment, the population $X$ is computed as the arithmetic average of the covering ratios across all hemispheres (**Section A.1**). The long-memory decay model parameterized by the gamma-type probability measure $F$ is then fitted to the experimental $X$ (**Table 2** and **Figure 4**), demonstrating a reasonable agreement between the experimental and theoretical results. **Table 2** indicates that the value of $\alpha$ is 0.2 to 0.3 and that of $\beta$ is approximately 1.

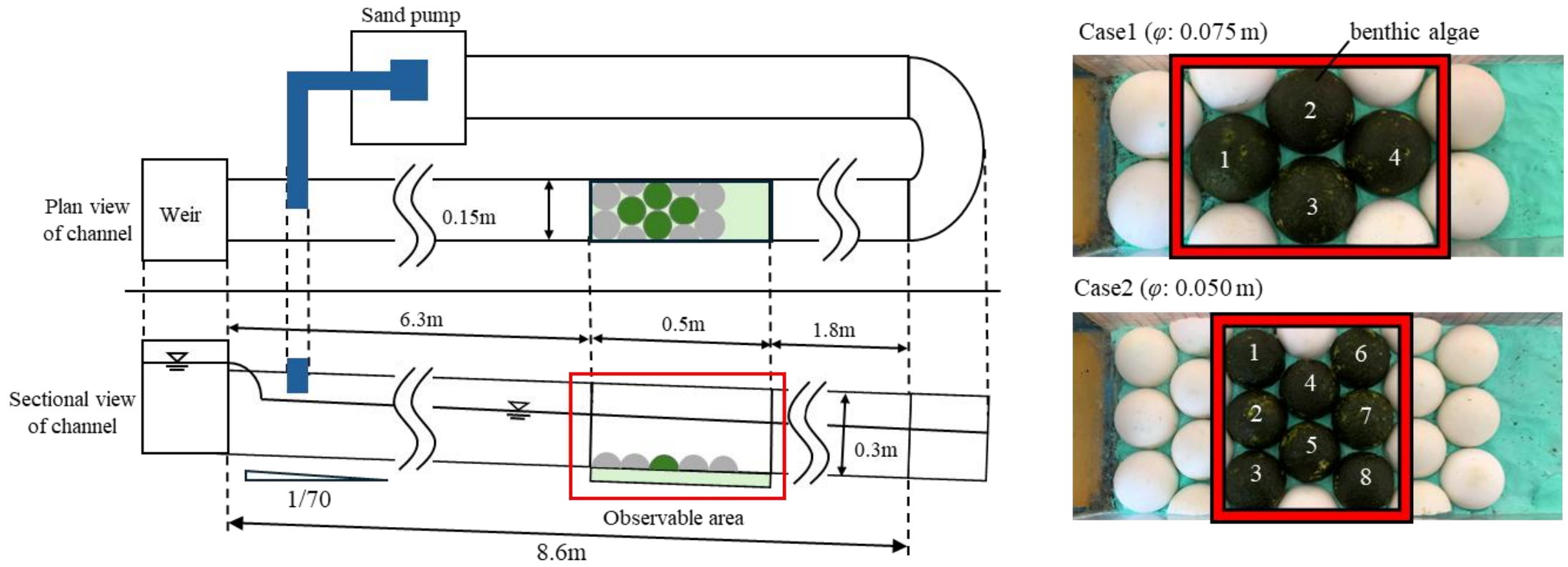

**Figure 3.** An image of the experimental setting. The left figure panel shows the experimental channel. The right panel shows the experimental configurations in the observable areas.

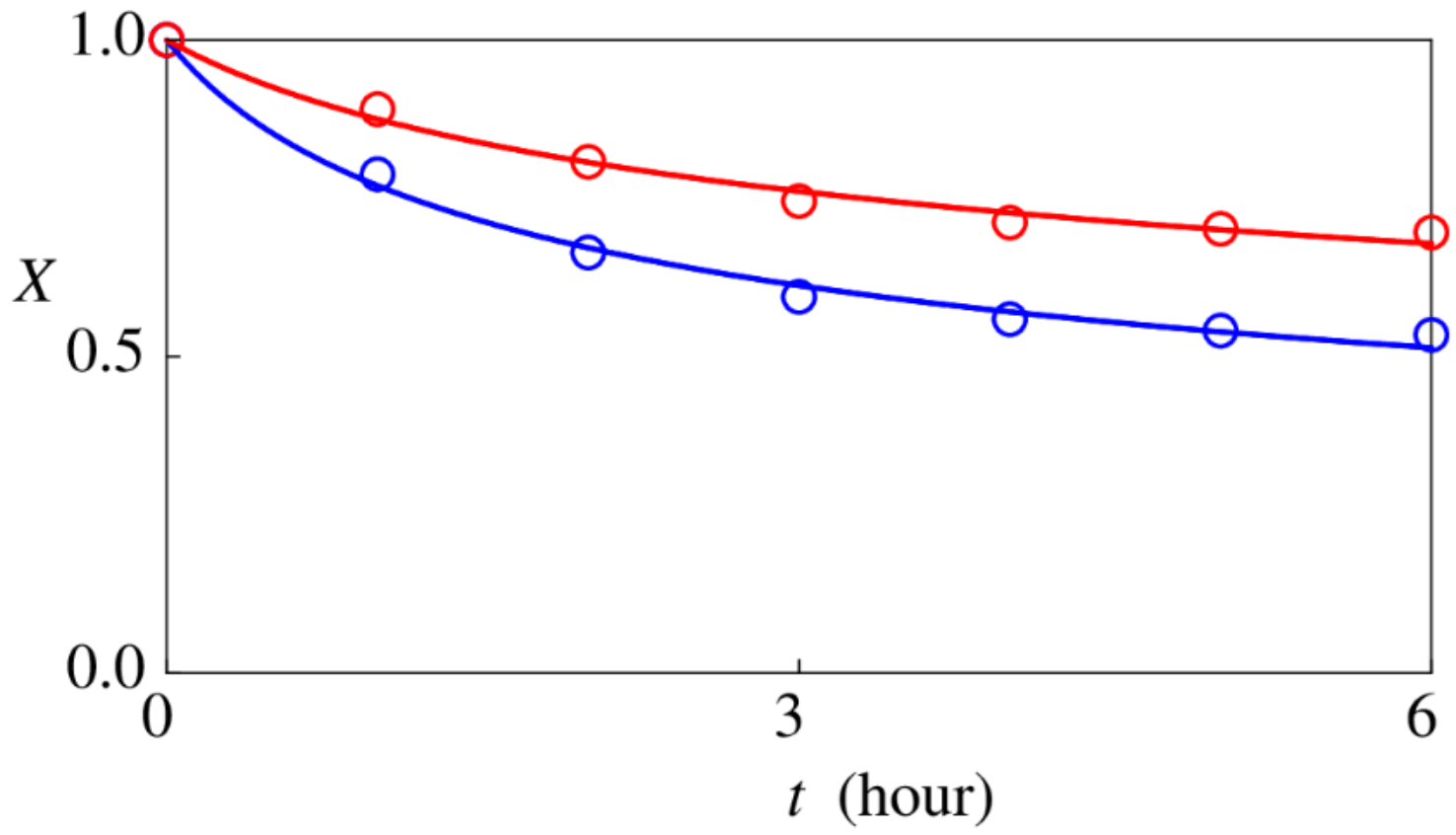

**Figure 4.** The decay of the benthic algae population is measured by the algae covering ratio of the surface of a hemisphere with the radius submerged in sediment-laden water flow. Circles denote data, Curves correspond to long-memory fit. Colors represent case 1 (blue) and case 2 (red). There is a good agreement between the theoretical and experimental results.

**Table 2.** Experimental setting and fitted parameter values.

| | case 1 | case 2 |
|---|---|---|
| Chanel width (m) | 0.15 | |
| Channel slope (m) | 1/70 | |
| Diameter of hemispheres (m) | 0.075 | 0.050 |
| Water discharge ($m^3/s$) | 0.0080 | 0.0065 |
| Diameter of sand particles | 0.002 | |
| Sand discharge ($m^3/s$) | 0.0000400 | 0.0000291 |
| $\alpha$ (-) | 0.2946 | 0.2103 |
| $\beta$ (1/hour) | 1.431 | 0.8881 |

### 3.2 Decay and growth

The computational analysis compares the microscopic stochastic system (19) and the macroscopic population dynamic model (21) to elucidate their convergence. The two models are computed with a common time increment $\Delta t$ using fully explicit discretization methods.

#### 3.2.1 Model setting

Consider the growth dynamics with

$$G(x) = rx(1-x)(x-a), \quad x \geq 0 \tag{30}$$

with a constant $a \in (0,1)$. This $G$ can be decomposed as follows (regularization is not introduced here):

$$g^{+}(x) = x \quad \text{and} \quad g^{-}(x) = a \,. \tag{31}$$

An elementary stability analysis indicates that the corresponding ODE (11) has three stationary solutions: stable equilibria at $x = 0,1$ and a saddle point at $x = a$. The incorporation of long-memory decay into this ODE eliminates or decreases the stable equilibrium $x = 1$, as determined by the right-hand side of (21), while preserving the equilibrium at $x = 0$. In the present setting, the decay due to abrasion would dominate the growth during high-flow events such as floods. In contrast, growth would dominate during low-flow events with a rich nutrient supply. The growth dynamics can be understood due to nutrient supply and competition among the population, while abrasion is an external physical factor acting on the population. Their equilibria therefore arise due to biological, ecological, and physical phenomena.

A positive equilibrium of (21) for $t = +\infty$ in this case, if it exists, must satisfy

$$\hat{x}_{\infty}(R) = \frac{r\hat{X}_{\infty} g^{+}\left(\hat{X}_{\infty}\right)}{R + r\hat{X}_{\infty} g^{+}\left(\hat{X}_{\infty}\right) + r\left(1-\hat{X}_{\infty}\right) g^{-}\left(\hat{X}_{\infty}\right)} = \frac{\hat{X}_{\infty}^{2}}{Rr^{-1} + \hat{X}_{\infty}^{2} - a\hat{X}_{\infty} + a}, \quad R > 0\,, \tag{32}$$

yielding the consistency equation to determine $\hat{X}_{\infty} > 0$:

$$1 = \int_{0}^{+\infty} \frac{\hat{X}_{\infty}}{Rr^{-1} + \hat{X}_{\infty}^{2} - a\hat{X}_{\infty} + a} F(\mathrm{d}R) \; \left(= H\left(\hat{X}_{\infty}\right)\right) \tag{33}$$

with

$$\frac{\mathrm{d}H\left(\hat{X}_{\infty}\right)}{\mathrm{d}\hat{X}_{\infty}} = \int_{0}^{+\infty} \frac{Rr^{-1} + a - \hat{X}_{\infty}^{2}}{\left(Rr^{-1} + \hat{X}_{\infty}^{2} - a\hat{X}_{\infty} + a\right)^{2}} F(\mathrm{d}R)\,. \tag{34}$$

It follows that $H(0) = 0$ and $H(1) < 1$, and $H$ admits exactly one maximum value in $(0,1)$ that is larger than 1 when $a$ and the average of $R$ are sufficiently small. In this case, (33) admits two solutions, and the larger one denoted by $\hat{X}_{\infty} = S_{\infty}$ corresponds to a stable equilibrium and the smaller one to a saddle. The corresponding $\hat{x}_{\infty}$ is obtained by substituting $\hat{X}_{\infty} = S_{\infty}$ into (32). When $a$ or the average of $R$ is sufficiently large, (33) has no solution in $(0,1)$ and only the state $\hat{X}_{\infty} = \hat{x}_{\infty} \equiv 0$ becomes the equilibrium.

Since our computational experiments aim to compare the stochastic system and population dynamics, particularly the convergence of the former to the latter as $M \in \mathbb{N}$ increases, the following parameter values are used unless otherwise specified: $r = 0.3/24$ (1/hour) and $a = 0.25$, where the growth rate is set considering the literature [14] and the reference therein (see Figure 1 and Table 1 in this literature), typically ranging between 0.2 and 0.3 (1/day). The value of $a$ here is hypothetical. The gamma-type $F$ is assumed with parameter values of case 1 presented in **Table 2**. The computational results exhibit qualitative consistency when the parameter values of case 2 are employed.

### 3.2.2 Numerical discretization

For a given $M$, both the stochastic system and population dynamics model are discretized in time using a classical forward Euler method with a time increment $\Delta t > 0$. The stochastic system (19) is discretized for each $i = 1,2,3,...,M$ and $k = 0,1,2,...$ as

$$x_{i,(k+1)\Delta t} = x_{i,k\Delta t} + \mathbb{I}\left(x_{i,k\Delta t} = 0\right)\Delta N_{i,k\Delta t} - \mathbb{I}\left(x_{i,k\Delta t} = 1\right)\Delta L_{i,k\Delta t} - \mathbb{I}\left(x_{i,k\Delta t} = 1\right)\Delta K_{i,k\Delta t} \tag{35}$$

with an initial condition $x_{i,0} \in \{0,1\}$. Here, $\Delta N_{i,k\Delta t}$ is a point process with the increment 1 and the jump rate of $X_{k\Delta t}^{(M)} g^{+}\left(X_{k\Delta t}^{(M)}\right)$, where $X_{k\Delta t}^{(M)} = \frac{1}{M}\sum_{i=1}^{M} x_{i,k\Delta t}$. Similarly, $\Delta L_{i,k\Delta t}$ is a point process with the increment 1 and the jump rate of $\left(1 - X_{k\Delta t}^{(M)}\right) g^{-}\left(X_{k\Delta t}^{(M)}\right)$, and $\Delta K_{i,k\Delta t}$ is a point process with the increment 1 and the jump rate of $R_i$. Each of $\Delta N_{i,k\Delta t}$, $\Delta L_{i,k\Delta t}$, and $\Delta K_{i,k\Delta t}$ are assumed mutually independent.

The population dynamics model (21) is discretized for $i = 1,2,3,...,M$ and $k = 0,1,2,...$ as

$$\hat{x}_{i,(k+1)\Delta t} = \hat{x}_{i,k\Delta t} + \Delta t\left[\hat{X}_{k\Delta t}^{(M)} g^{+}\left(\hat{X}_{k\Delta t}^{(M)}\right) - \left\{R_i + \hat{X}_{k\Delta t}^{(M)} g^{+}\left(\hat{X}_{k\Delta t}^{(M)}\right) + \left(1 - \hat{X}_{k\Delta t}^{(M)}\right) g^{-}\left(\hat{X}_{k\Delta t}^{(M)}\right)\right\}\hat{x}_{i,k\Delta t}\right] \tag{36}$$

starting from an initial condition $\hat{x}_{i,0} \in [0,1]$, where $\hat{X}_{k\Delta t}^{(M)} = \frac{1}{M}\sum_{i=1}^{M} \hat{x}_{i,k\Delta t}$. Here, $\hat{x}_{i,k\Delta t}$ is an approximation of $\hat{x}_t\left(R_i\right)$ at time $k\Delta t$. The time increment is $\Delta t = 0.001$ (day) that has been found to be sufficiently fine for our computational purpose. The initial condition is fixed to be 1 for both the stochastic system and population dynamics model. Below, the hat for $x$ and $X$ are omitted for simplicity.

### 3.2.3 Results and discussion

Convergence behavior of the stochastic system to the population dynamics model (**Figure 5**) is investigated. **Figure 5(a)** visualizes the population $X$ computed from the stochastic system with $M = 2^l$ ($l = 1,2,3,...,16$) where the growth dynamics are neglected by setting $r = 0$. Similarly, **Figure 5(b)** visualizes the results that account for the growth dynamics. As depicted in **Figure 5(a)**, all the computed trajectories are decreasing in time when the growth dynamics are not considered. By contrast, in **Figure 5(b)**, the computed trajectories of $X$ are non-monotone in time due to the presence of the growth dynamics that emerge as positive jumps. Positive jumps in the trajectories of the stochastic system become less significant as $M$ increases.

The convergence speed from the stochastic system to the population dynamics model is quantified by computing the average of the error $\left(X_{t,\text{micro}} - X_{t,\text{macro}}\right)^2$ at $t = k\Delta t$ ( $k = 1, 2, 3, ..., 7000$ ), where $X_{t,\text{micro}}$ and $X_{t,\text{macro}}$ are the populations by the stochastic system and population dynamics model, respectively. The least-squares error $\text{Er} = \text{Er}\left(M\right)$ as a function of $M$ is plotted in **Figure 6**. The linear regressed curves of $\text{Er}\left(M\right)$ for $M = 2^l$ by a least-squares between the computed and fitted curves (on a common logarithmic scale for $\text{Er}$ ) are $0.092 \times 2^{-1.06l}$ ( $\text{R}^2$ value is 0.851) for the case without growth and $0.1023 \times 2^{-1.02l}$ ( $\text{R}^2$ value is 0.882) for that with growth, respectively. The fitted results thus suggest that the deviation of the population $X$ between the stochastic system and population dynamics model decay as $O\left(M^{-1/2}\right)$. The obtained computational results support the convergence of the stochastic system to the population dynamics model as $M$ increases. Moreover, this dependence on $M$ is in accordance with the theoretical estimate (49).

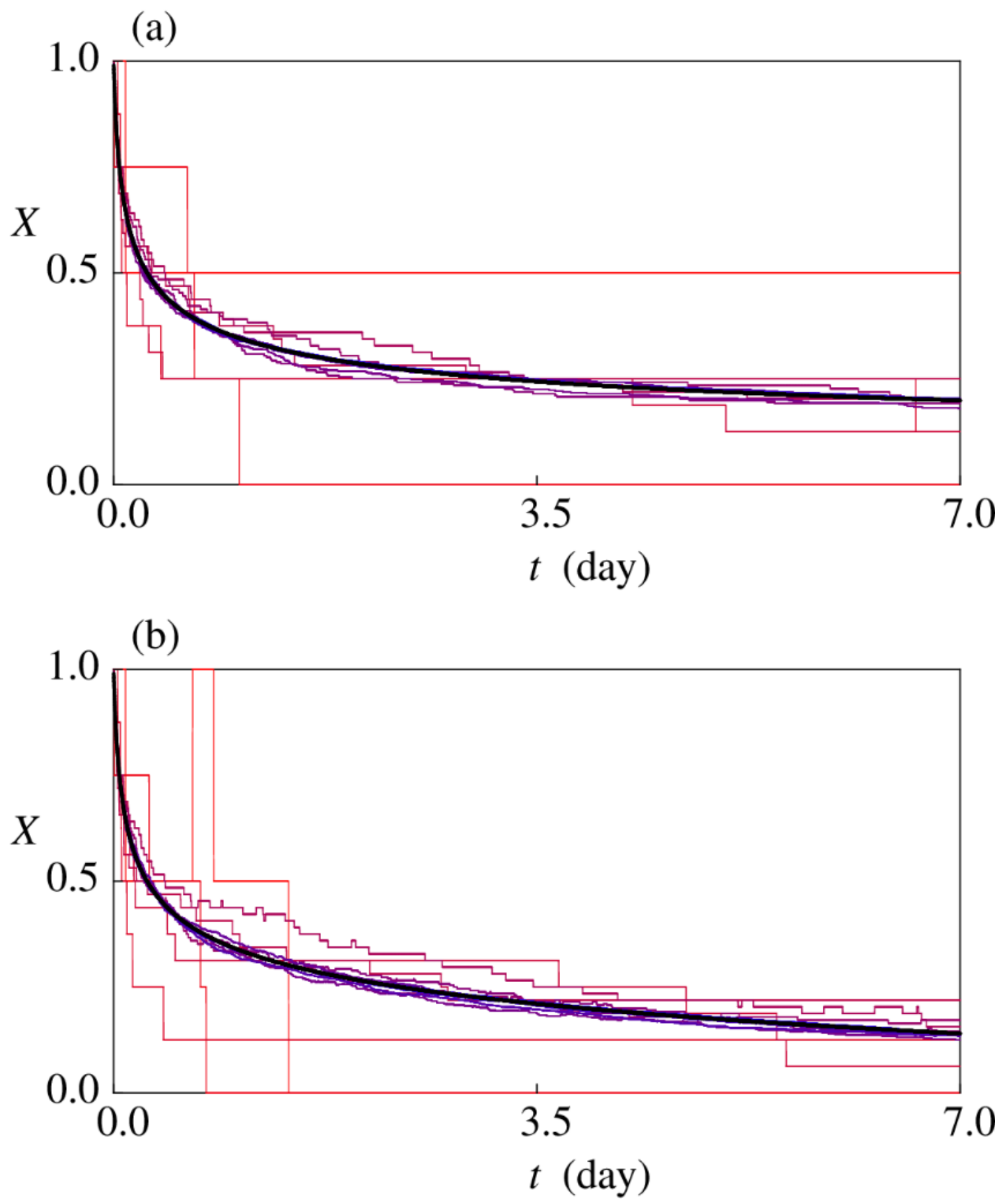

**Figure 5.** Convergence visualization of the population $X$ from the stochastic system (colored from red to blue as $M$ increases from 0 to $2^{16}$) to the population dynamics model (black curve): (a) a case with decay but without growth and (b) a case with both growth and decay. Both figure panels show the convergence behavior of spin processes.

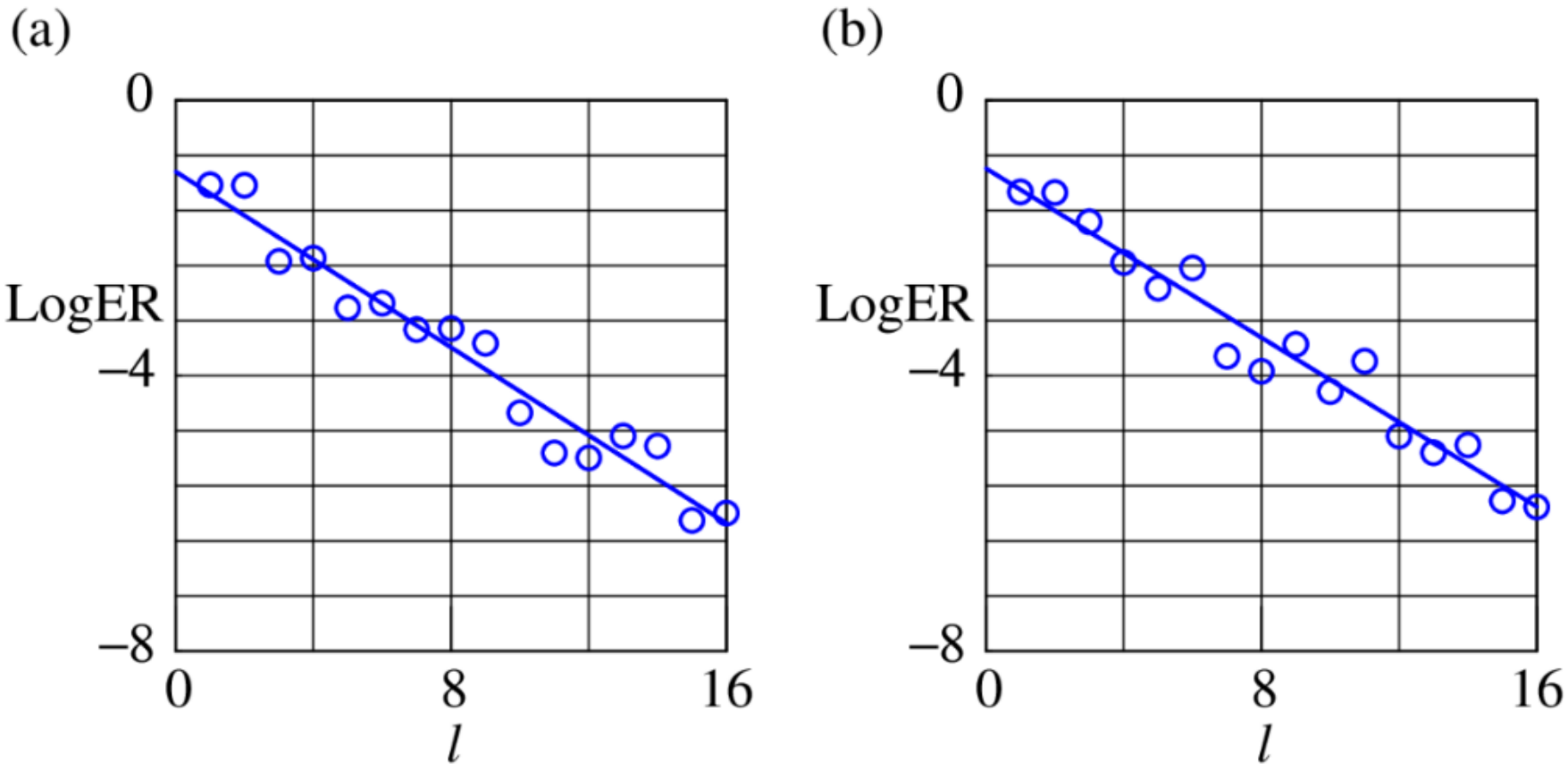

**Figure 6.** Logarithm of the least-squares errors (ER) of the population $X$ between the stochastic system to the population dynamics model based on the computational results (circles) and least-squares fitting (curves): (a) a case with decay without growth and (b) a case with both growth and decay. Both figure panels show that the regression serves well.

### 3.3 Application to rate-induced tipping

As an extended version of the proposed mathematical models, this paper computationally analyzes the growth dynamics associated with $G$ of (30), where the coefficient $a$ is a positive, time-dependent parameter. In these models, a temporally varying $a$ induces temporal shifts in the saddle point separating two stable equilibria (if they exist). From an engineering perspective, this corresponds to a situation where decreasing $a$ (resp., increasing $a$) signifies some improvement (resp., degradation) of the conditions for providing the algae population, such as nutrient availability and light exposure [73,74,75]. This subsection emphasizes interactions between population decay caused by abrasion and biological growth influenced by dynamic environmental changes. The computation in this section is exploratory but would give insights into more advanced models in the future.

In the sequel, the coefficient $a_t$, which now depends on time $t \geq 0$, is assumed:

$$a_t = \bar{a} + \frac{\underline{a} - \bar{a}}{2}\left\{1 + \tanh\left(\frac{t-h}{\theta}\right)\right\}, \tag{37}$$

where $\underline{a}, \bar{a}$ with $0 < \underline{a} < \bar{a} < 1$ are the lower and upper bounds of $a_t$, $h > 0$ is the shift, and $\theta > 0$ is the scaling. This $a_t$ is a decreasing sigmoidal curve connecting $\underline{a}$ and $\bar{a}$. The parameter values are set as $\underline{a} = 0.1$, $\bar{a} = 0.5$, $h = 30$ (day), and $\theta = 2$ (day). With these parameterizations, the time-dependent $a$ has an inflection point at $t = h$ and has a profile such that it is close to $\bar{a}$ for $t$ sufficiently small, and is close to $\underline{a}$ for $t$ sufficiently large (**Figure 7**). These parameter values are hypothetical, while they can cover cases with and without tipping phenomenon.

Sigmoidal curves of the form (37) serve as a representative model for temporal environmental fluctuations, facilitating the efficient analysis of rate-induced tipping [76,77]. For ODEs, the rate-induced tipping in this context arises when model parameters satisfy some threshold condition, resulting in population $X$ extinction ($X_t \to 0$ as $t \to +\infty$) or persistence ($X_t \to S_\infty\big|_{a=\underline{a}}$ as $t \to +\infty$) depending on whether the condition is violated. Specifically, considering $X_0 = 1$, if there is no population decay (if $X_t > a_t$ for all $t > 0$) then there is no tipping, while the tipping occurs if $X_t < a_t$ at some $t > 0$ (e.g., Feudel [53]); the latter occurs for example if $r$ is sufficiently small. However, the proposed model introduces a more complex tipping mechanism as the population is an aggregated variable being different from the case of classical ODEs. Additionally, the tipping dynamics of the stochastic system, influenced by the degree of freedom $M$, are critical for understanding the effect of $M$ on system stability.

What is (indirectly) controllable in the present setting would be the distribution $F$ by changing the scheme of sediment supply, e.g., river discharge and sediment diameter. Changing the abrasion would affect the decay dynamics of the benthic algae population and hence their tipping behavior. First, the system stability for different values of $\beta$ in the probability measure $F$ is investigated, where the replacement $\beta \to \eta\beta$ is employed in $F$ to represent the situation where a larger $\beta$ implies more rapid abrasion of the benthic algae, and vice versa. The nominal case corresponds to $\eta = 1$. A particular interest in applications would be what will happen if the abrasion becomes weaker or the flow velocity decreases due to lowering the discharge or channel slope by decreasing the sediment discharge, namely, by decreasing $\eta$. Such quantitative experimental data are currently not available, and therefore, the computational experiments presented below only give theoretical insights; nevertheless, they provide the first case study of rate-induced tipping of the mathematical models of the proposed forms.

**Figure 8** illustrates the computed trajectories of the population dynamics model with $\eta = 1$, $0.0094$, $0.0093$, and the time-dependent coefficient $a$. According to **Figure 8**, there exists a threshold value of $\eta = \eta_c$ between 0.0093 and 0.0094, such that the population goes extinct if $\eta > \eta_c$ owing to the sufficiently strong abrasion dominating population growth, while the population eventually approaches a positive value if $\eta < \eta_c$ as growth dominates over decay. Notably, the positive equilibrium is eventually attainable even when the population $X_t$ crosses $a_t$, indicating that the classical theory of rate-induced tipping for ODEs does not apply to our model. These computational results emphasize that the balance between decay (or abrasion) and growth should be carefully evaluated in field applications, as small modeling errors may lead to significantly different long-term predictions for benthic algae populations. More specifically, from an application perspective, sediment replenishment projects have been implemented in many rivers to artificially supply sediment particles to mitigate sediment depletion caused by dam/weir construction or mining [78,79]. In this context, supplying a larger amount of sediments corresponds to assigning a larger value of $\eta$. Cost-efficient sediment supply aimed at suppressing the bloom of nuisance benthic algae involves identifying the critical value $\eta = \eta_c$. The value of $\eta$ is determined on a case-specific basis; however, the population dynamics model facilitates the determination of its value through computational experiments, as demonstrated in this study.

Finally, the rate-induced tipping of the stochastic system to analyze the influence of stochasticity on system stability is investigated. The parameter values are $M = 128$, and the stochastic system is computed with the coefficient $a_t$ specified in (37). Several values of $\eta$ from approximately 0.0093 to 0.0094 are then examined, as the critical value $\eta = \eta_c$ for the population dynamics model lies within this interval. To computationally examine the system stability, the stochastic system over a sufficiently long duration, which is 200 (day), and the histogram of the population $X$ at this time are computed. The shape of the histogram, such as the total number of maximum points, reflects the stability of the stochastic system. For each value of $\eta$, 10,000 sample paths of the stochastic system are generated.

**Figure 9** illustrates the computed histograms corresponding to different values of $\eta$ near the critical value $\eta = \eta_c$ for the population dynamics model. For each value of $\eta$ considered, the histogram displays two modes: one at zero indicating population extinction, and another between 0.75 and 0.80. The latter diminishes significantly for $\eta = 0.020$. According to **Figure 9**, the existence of stochasticity mitigates the sharp tipping phenomenon observed in the deterministic case, and the attractivity of the two stable equilibria is depicted by the height of the histograms. Furthermore, the stability of the equilibrium associated with extinction increases with abrasion, i.e., $\eta$.

To better comprehend the influences of the degree of freedom $M$, the histograms against different values of $M$ with the fixed value $\eta = 0.008$ are computed, as illustrated in **Figure 10**. The computed histograms become sharper as $M$ increases and as the stochastic system approaches the population dynamics model. The bimodality of the histograms remains invariant within the range $M = 2^7$ to $M = 2^{10}$; however, the peak at $X = 0$ diminishes as $M$ increases. This peak is anticipated to vanish in the limit $M = +\infty$, as suggested in **Figure 9**, due to $0.008 < \eta_c$. Consequently, the tipping behavior of the stochastic system is expected to approximate that of the population dynamics model in a weak sense such that the histograms in **Figure 10** converge to Dirac's delta concentrated at a positive point.

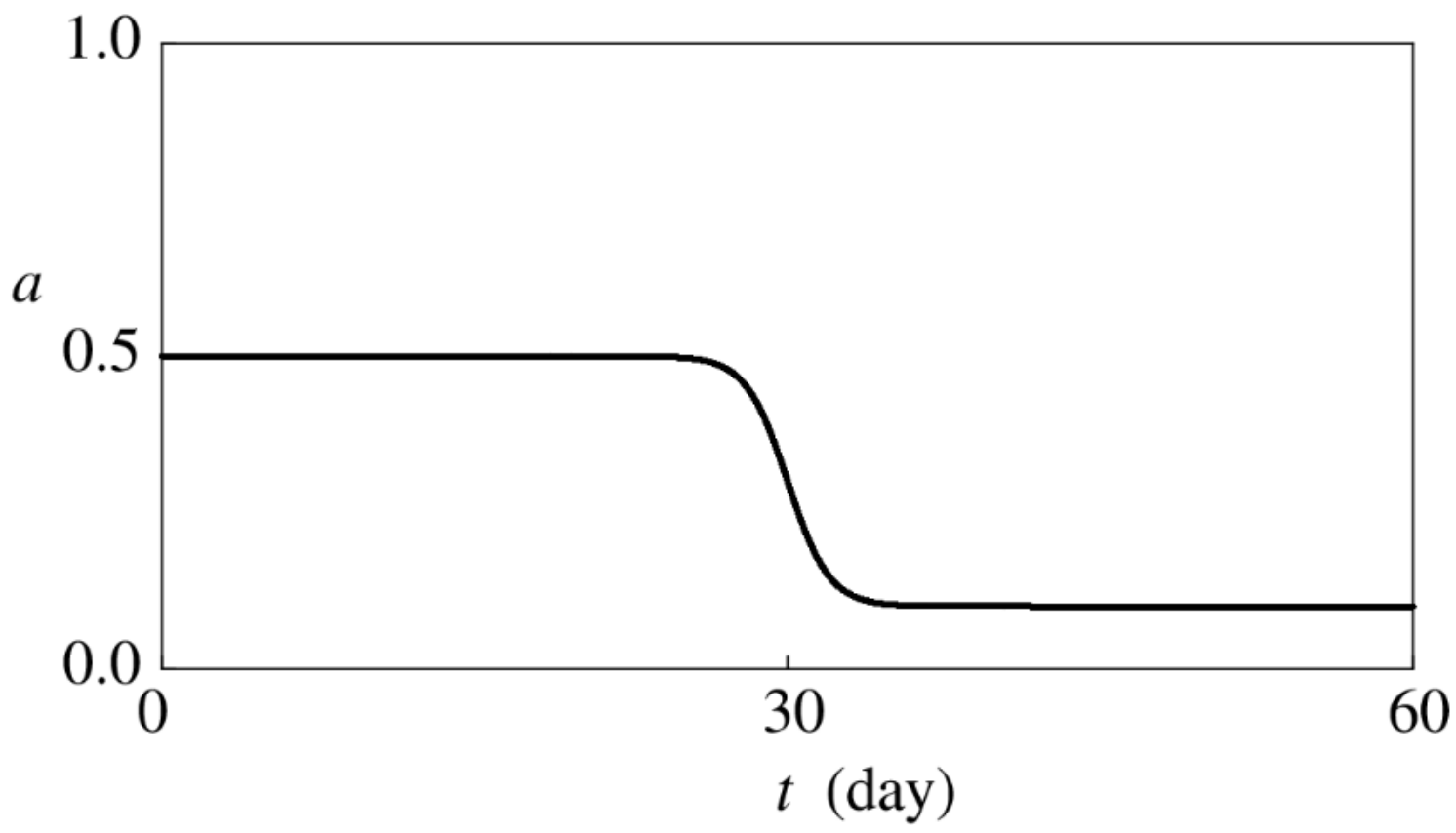

**Figure 7.** The specified profile of coefficient $a$, which has a decreasing sigmoidal shape.

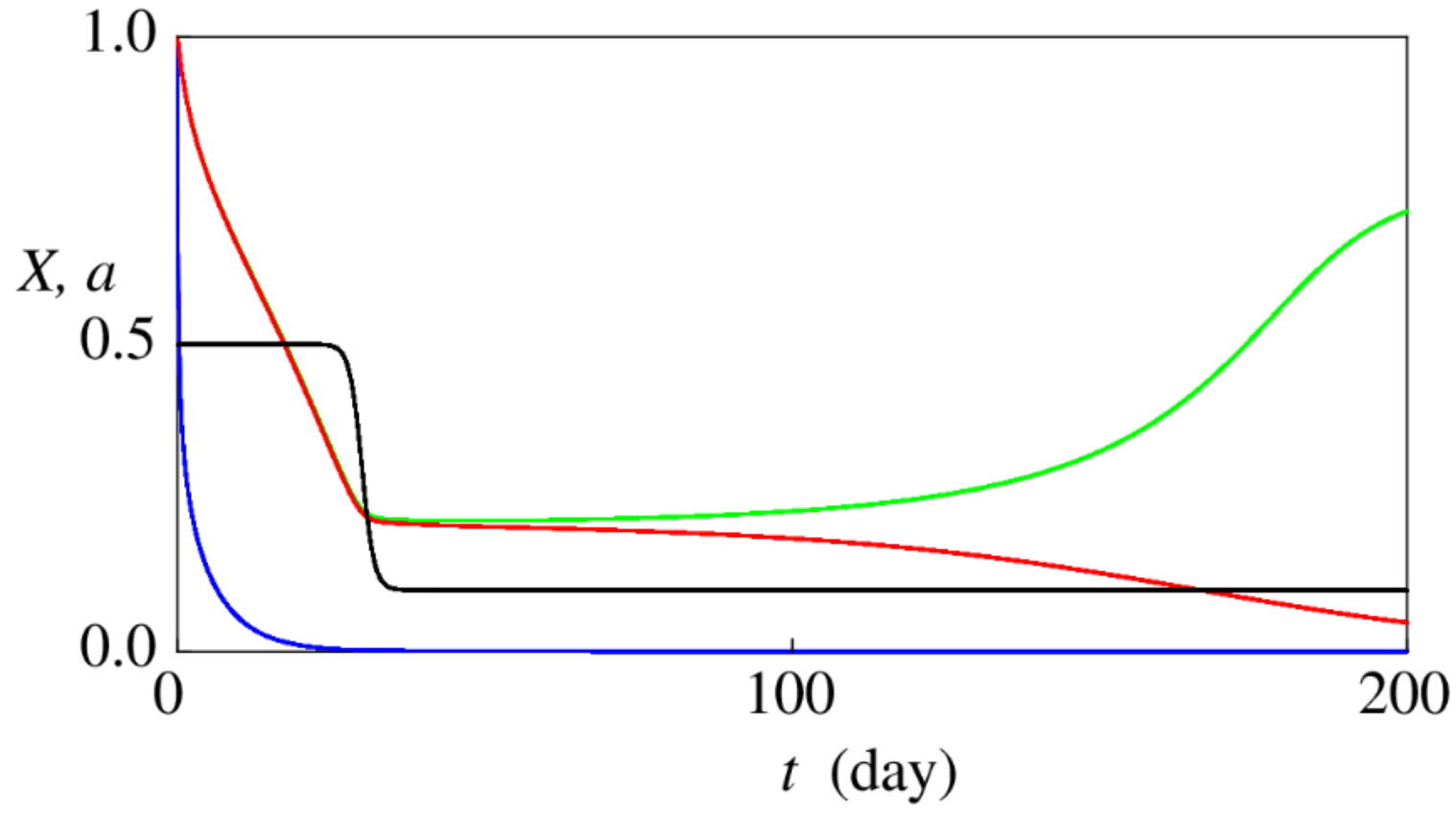

**Figure 8.** Population profiles $X$ modeled for $\eta = 1$ (blue, the nominal case), $\eta = 0.0094$ (red), $\eta = 0.0093$ (green), and the time-dependent sigmoidal coefficient $a$ (black). The green curve corresponds to the scenario converging to the positive stable equilibrium $S_\infty$ while the red one to the stable zero equilibrium. The occurrence of tipping is visible by comparing the red and green curves.

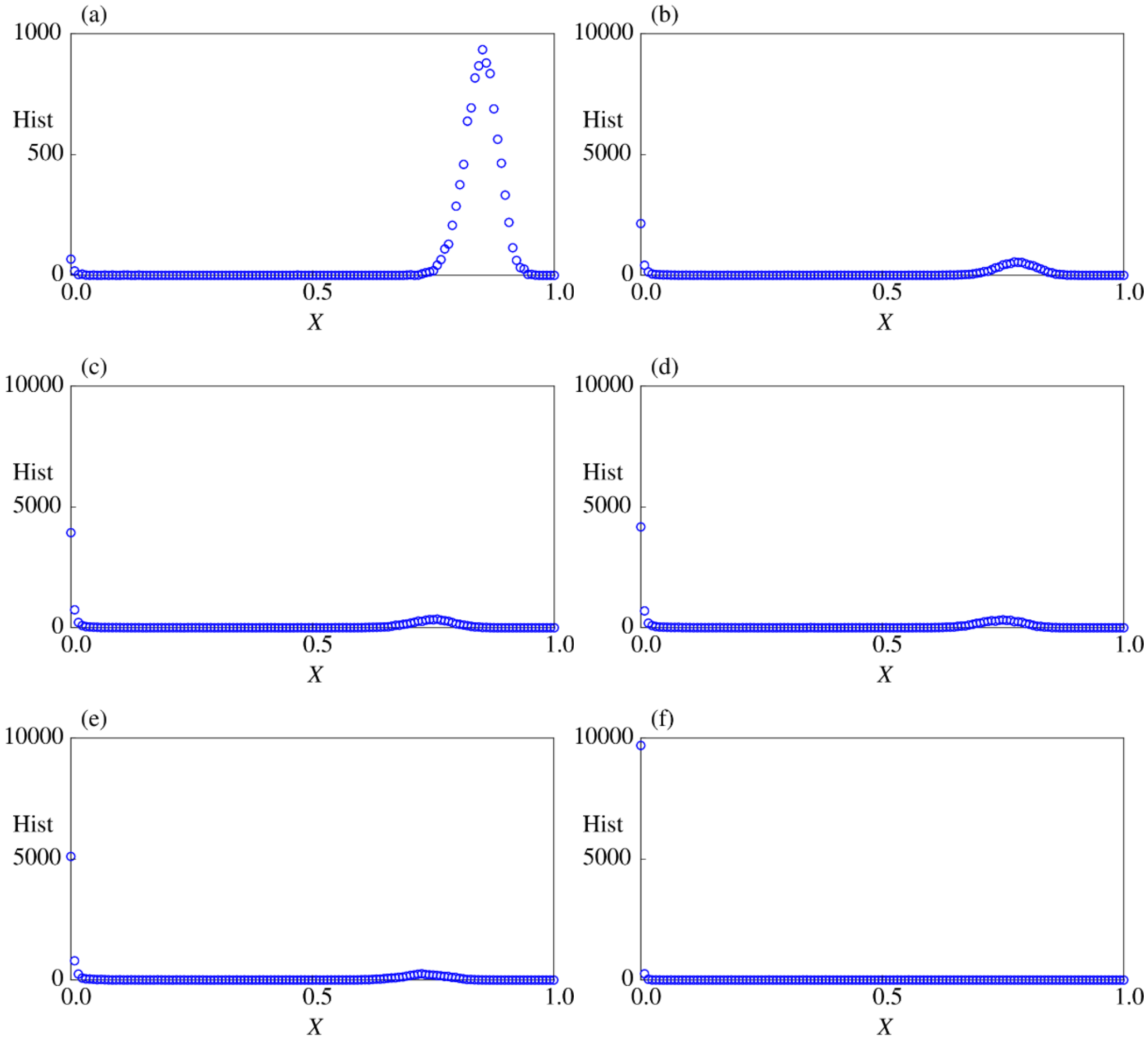

**Figure 9.** Computed histograms (Hist) for various values of $\eta$, based on 10,000 sample paths: (a) $\eta = 0.005$, (b) $\eta = 0.008$, (c) $\eta = 0.0093$, (d) $\eta = 0.0094$, (e) $\eta = 0.010$, and (f) $\eta = 0.020$. Increasing the strength of abrasion ($\eta$) leads to a profile more concentrated at the origin.

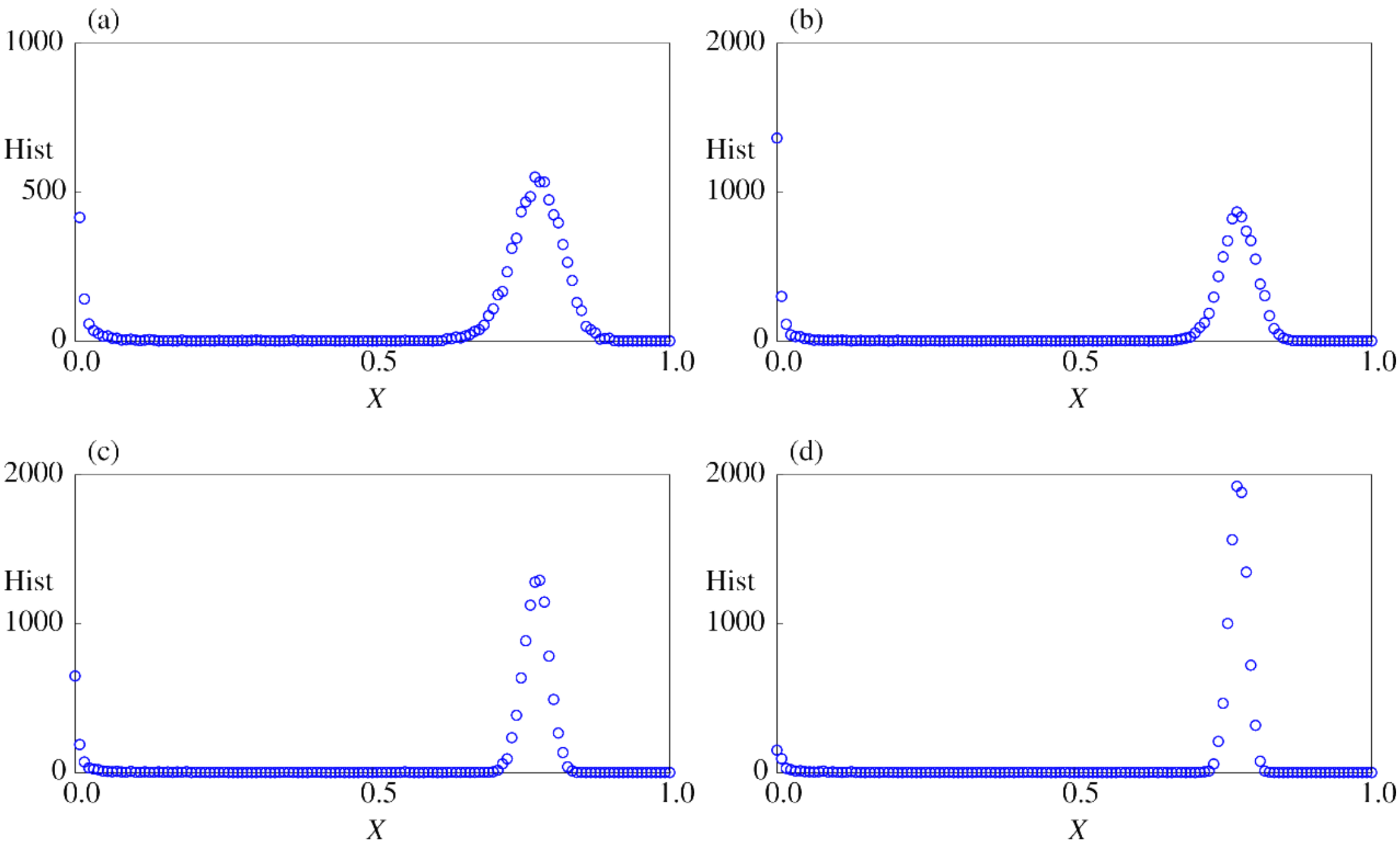

**Figure 10.** Computed histograms (Hist) for various values of $M$, based on 10,000 sample paths: (a) $M = 2^7$, (b) $M = 2^8$, (c) $M = 2^9$, and (d) $M = 2^{10}$. Increasing the total number of sample paths yields a more concentrated distribution.

## 4. Conclusions

This paper proposed a unified mathematical framework for describing both the growth and decay of benthic algae population dynamics. The approach employs appropriate superpositions of microscopic heterogeneous spin processes to derive the macroscopic population dynamics. The resulting governing equation characterizes the population dynamics as a continuum of ODEs, being distinct from classical ODEs. Computational examples were presented to illustrate the proposed model and its applications to rate-induced tipping phenomena, although the latter was addressed computationally. Theoretically, a wide class of growth curves can be modeled within the proposed mathematical framework, as demonstrated in **Section 2.3**. Moreover, a model with an exponential decay can be derived if is $F$ set to be a suitable Dirac delta.

Our theoretical framework can be extended to the case where the population dynamics are driven by an external noise process representing transient river flows. It has been demonstrated that streamflow dynamics are non-Markovian and exhibit long-memory behavior [80,81]. Thus, a more realistic model of the population decay of benthic algae would involve long-memory decay driven by another long-memory process, which will result in a nontrivial population dynamics model. Indeed, a limitation of our study is the lack of mathematical linkages between physical and model parameters due mainly to the conceptual nature of the proposed model. Establishing a physical connection between the probability measure $F$ of spatial heterogeneity and the subdomain locations of $D$, would enable the proposed model to address a broader range of riverbed geometries. These aspects are currently under investigation by the authors so that ecological and biological implications of the proposed model can be enhanced. The tipping phenomenon will be studied theoretically in detail in future in view of dynamical systems having long memory. Employing a machine-learning method [82] would enhance applicability of the proposed model to real-world problems. Application of the proposed model to a field case study will also be an important step towards the development of an advanced modeling framework for the population dynamics of benthic algae. An interesting application in this direction will be the development for coupled nutrient-algae dynamics [83].

Conducting field surveys in a short time is difficult, and long-term data acquisition is required. On the other hand, as mentioned above, the experimental fact that the decay rate of attached algae is not exponential is expected to be critical for civil engineering applications, and the authors believe that the approach presented in this paper may be a starting point for overcoming the limitations of previous research. Of course, this study has limitations, such as not using field data. However, even considering the issue of physical scale, it is thought that the question of whether the decay rate is exponential or not will not qualitatively change, suggesting a certain degree of engineering applicability of the proposed framework. While the uniqueness of this perspective and approach at present is both an advantage and a disadvantage, its true value is expected to become even clearer as more field data is accumulated in the future.

## Appendix

### A.1 Experimental data

For each experiment (case 1 and case 2 in the main text), this section presents the observed population $X_t$ at each discrete time step $t$ as the arithmetic average of the covering ratios across all hemispheres, as depicted in **Tables A1** and **A2**. See also **Figure 4** in the main text.

**Table A1**. The observed population $X_t$ at each discrete time step $t$ of each hemisphere for case 1.

| | | Hemisphere | | | |
|---|---|---|---|---|---|
| Time (s) | Average | 1 | 2 | 3 | 4 |
| 0 | 1.00.E+00 | 1.00.E+00 | 1.00.E+00 | 1.00.E+00 | 1.00.E+00 |
| 3600 | 7.88.E-01 | 9.68.E-01 | 7.97.E-01 | 6.84.E-01 | 7.02.E-01 |
| 7200 | 6.64.E-01 | 9.03.E-01 | 6.52.E-01 | 5.48.E-01 | 5.53.E-01 |
| 10800 | 5.94.E-01 | 7.93.E-01 | 5.90.E-01 | 5.05.E-01 | 4.89.E-01 |
| 14400 | 5.59.E-01 | 7.46.E-01 | 5.53.E-01 | 4.77.E-01 | 4.58.E-01 |
| 18000 | 5.41.E-01 | 7.16.E-01 | 5.40.E-01 | 4.67.E-01 | 4.38.E-01 |
| 21600 | 5.34.E-01 | 7.05.E-01 | 5.31.E-01 | 4.67.E-01 | 4.32.E-01 |

**Table A2**. The same is true for **Table A1** for case 2.

| | | Hemisphere | | | | | | | |
|---|---|---|---|---|---|---|---|---|---|
| Time (s) | Average | 1 | 2 | 3 | 4 | 5 | 6 | 7 | 8 |
| 0 | 1.00.E+00 | 1.00.E+00 | 1.00.E+00 | 1.00.E+00 | 1.00.E+00 | 1.00.E+00 | 1.00.E+00 | 1.00.E+00 | 1.00.E+00 |
| 3600 | 8.91.E-01 | 9.74.E-01 | 9.49.E-01 | 9.72.E-01 | 9.28.E-01 | 9.65.E-01 | 9.12.E-01 | 9.03.E-01 | 5.22.E-01 |
| 7200 | 8.07.E-01 | 9.25.E-01 | 9.03.E-01 | 8.82.E-01 | 7.92.E-01 | 8.51.E-01 | 8.04.E-01 | 8.43.E-01 | 4.59.E-01 |
| 10800 | 7.45.E-01 | 8.69.E-01 | 8.38.E-01 | 8.31.E-01 | 7.38.E-01 | 7.77.E-01 | 7.08.E-01 | 8.00.E-01 | 4.02.E-01 |
| 14400 | 7.11.E-01 | 8.43.E-01 | 7.92.E-01 | 7.69.E-01 | 7.06.E-01 | 7.41.E-01 | 6.85.E-01 | 7.74.E-01 | 3.82.E-01 |
| 18000 | 7.01.E-01 | 8.38.E-01 | 7.85.E-01 | 7.53.E-01 | 6.91.E-01 | 7.29.E-01 | 6.75.E-01 | 7.69.E-01 | 3.71.E-01 |
| 21600 | 6.96.E-01 | 8.34.E-01 | 7.80.E-01 | 7.41.E-01 | 6.88.E-01 | 7.22.E-01 | 6.66.E-01 | 7.66.E-01 | 3.70.E-01 |

### A.2 Proofs of Propositions

***Proof of Proposition 1***

The proof proceeds by demonstrating that the average of the left-hand side of (8) equals the right-hand side, while the variance of the right-hand side of (8) vanishes. The technical aspects requiring further classification are detailed below. Time $t \geq 0$ is fixed in the sequel.

The average and variance of $X_t^{(M)}$ are calculated directly using (7) as

$$\mathbb{E}\left[X_t^{(M)}\right] = \frac{1}{M}\sum_{i=1}^{M}\mathbb{E}\left[x_{i,t}\right] = \frac{1}{M}\sum_{i=1}^{M} e^{-R_i t} \tag{38}$$

and

$$\begin{aligned}\mathbb{V}\left[X_t^{(M)}\right] &= \mathbb{E}\left[\left(\frac{1}{M}\sum_{i=1}^{M} x_{i,t}\right)^2\right] - \left(\frac{1}{M}\sum_{i=1}^{M}\mathbb{E}\left[x_{i,t}\right]\right)^2 \\ &= \frac{1}{M^2}\mathbb{E}\left[\sum_{i,j=1}^{M} x_{i,t}x_{j,t}\right] - \frac{1}{M^2}\sum_{i=1}^{M} e^{-2R_i t} - \frac{2}{M^2}\sum_{i<j}^{M} e^{-R_i t}e^{-R_j t}\end{aligned} . \tag{39}$$

By $\left(x_{i,t}\right)^2 = x_{i,t}$ and independence between $x_{i,t}, x_{j,t}$ ($i \neq j$),

$$
\begin{aligned}
\mathbb{E}\left[\sum_{i,j=1}^{M} x_{i,t}x_{j,t}\right] &= \mathbb{E}\left[\sum_{i=1}^{M}\left(x_{i,t}\right)^2\right] + \mathbb{E}\left[2\sum_{i<j}^{M} x_{i,t}x_{j,t}\right] \\
&= \sum_{i=1}^{M}\mathbb{E}\left[x_{i,t}\right] + 2\sum_{i<j}^{M}\mathbb{E}\left[x_{i,t}\right]\mathbb{E}\left[x_{j,t}\right] \\
&= \sum_{i=1}^{M} e^{-R_i t} + 2\sum_{i<j}^{M} e^{-R_i t}e^{-R_j t}
\end{aligned}
\tag{40}
$$

and hence

$$
\mathbb{V}\left[X_t^{(M)}\right] = \frac{1}{M^2}\sum_{i=1}^{M} e^{-R_i t}\left(1-e^{-R_i t}\right). \tag{41}
$$

Second, it is shown that the limit $\lim_{M\to+\infty} X_t^{(M)}$ exists and is not random, i.e., it is a constant. Due to (41) and

$$
0 \le \lim_{M\to+\infty}\frac{1}{M^2}\sum_{i=1}^{M} e^{-R_i t}\left(1-e^{-R_i t}\right) \le \lim_{M\to+\infty}\frac{1}{M^2}\sum_{i=1}^{M} 1 = 0, \tag{42}
$$

it follows that $\lim_{M\to+\infty}\mathbb{V}\left[X_t^{(M)}\right] = 0$. This shows that $X_t^{(M)}$ converges to a constant in the square norm and hence in the sense of probability due to the Markov inequality (e.g., Theorem 5.11 in Klenke [41]). This nonrandom limit is denoted by $\hat{X}_t$, which must be $\lim_{M\to+\infty}\mathbb{E}\left[X_t^{(M)}\right]$. Clearly, $X_t^{(M)}$ is bounded between 0 and 1. Then, the dominated convergence (Theorem 2.16 in Jiang [84]) shows

$$
\lim_{M\to+\infty}\mathbb{E}\left[X_t^{(M)}\right] = \mathbb{E}\left[\lim_{M\to+\infty} X_t^{(M)}\right] = \mathbb{E}\left[\hat{X}_t\right] = \hat{X}_t. \tag{43}
$$

Finally, by (38) and (43),

$$
\hat{X}_t = \lim_{M\to+\infty}\frac{1}{M}\sum_{i=1}^{M} e^{-R_i t} = \lim_{M\to+\infty}\int_0^{+\infty} e^{-Rt}F_M\left(\mathrm{d}R\right), \tag{44}
$$

and the proof concludes if the right-most side of (44) is determined. By virtue of the discretization of the domain $D$, it follows that (e.g., Proof of Proposition 2 in Yoshioka [62]):

$$
\left|F\left(\left(0,R\right)\right) - F_M\left(\left(0,R\right)\right)\right| \le \frac{1}{M} \quad \text{for all } R > 0, \tag{45}
$$

showing that $F_M$ uniformly approximates $F$. This implies that $F_M$ converges weakly to $F$ (Definition 13.21 and Theorem 13.23 in Klenke [41]) because $F$ is continuous due to **Assumption 1**. Then, by the Portemanteau theorem (Theorem 13.16(i)-(ii) in Klenke [41]), it follows that

$$
\lim_{M\to+\infty}\int_0^{+\infty} e^{-Rt}F_M\left(\mathrm{d}R\right) = \int_0^{+\infty} e^{-Rt}F\left(\mathrm{d}R\right), \tag{46}
$$

which combined with (44) yields (8).

□

***Proof of Proposition 2***

Fix $T > 0$. Then, fix $t \in \left[0,T\right]$. First, (16) is proven. A direct calculation shows

$$
\begin{aligned}
&\mathbb{E}\left[\left(\frac{1}{M}\sum_{i=1}^{M}\left\{\int_0^t\left(1-x_{i,s-}\right)\mathrm{d}\tilde{N}_{i,s}-\int_0^t x_{i,s-}\mathrm{d}\tilde{L}_{i,s}\right\}\right)^2\right]\\
&=\frac{1}{M^2}\sum_{i,j=1}^{M}\mathbb{E}\left[\begin{array}{l}\int_0^t\left(1-x_{i,s-}\right)\mathrm{d}\tilde{N}_{i,s}\int_0^t\left(1-x_{j,s-}\right)\mathrm{d}\tilde{N}_{j,s}+\int_0^t x_{i,s-}\mathrm{d}\tilde{L}_{i,s}\int_0^t x_{j,s-}\mathrm{d}\tilde{L}_{j,s}\\ -\int_0^t\left(1-x_{i,s-}\right)\mathrm{d}\tilde{N}_{i,s}\int_0^t x_{j,s-}\mathrm{d}\tilde{L}_{j,s}-\int_0^t x_{i,s-}\mathrm{d}\tilde{L}_{i,s}\int_0^t\left(1-x_{j,s-}\right)\mathrm{d}\tilde{N}_{j,s}\end{array}\right],\\
&=\frac{1}{M^2}\sum_{i,j=1}^{M}\mathbb{E}\left[\int_0^t\left(1-x_{i,s-}\right)\mathrm{d}\tilde{N}_{i,s}\int_0^t\left(1-x_{j,s-}\right)\mathrm{d}\tilde{N}_{j,s}+\int_0^t x_{i,s-}\mathrm{d}\tilde{L}_{i,s}\int_0^t x_{j,s-}\mathrm{d}\tilde{L}_{j,s}\right]\\
&=\frac{1}{M^2}\sum_{i=1}^{M}\mathbb{E}\left[\left(\int_0^t\left(1-x_{i,s-}\right)\mathrm{d}\tilde{N}_{i,s}\right)^2+\left(\int_0^t x_{i,s-}\mathrm{d}\tilde{L}_{i,s}\right)^2\right]
\end{aligned}
\tag{47}
$$

where the last equality follows from the independence between the Poisson random measures defining $\tilde{N}_i$ and $\tilde{N}_j$ and those defining $\tilde{L}_i$ and $\tilde{L}_j$ ($i \neq j$).

By considering the jump rates of each $N_i$ and $L_i$,

$$
\begin{aligned}
&\frac{1}{M^2}\sum_{i=1}^{M}\mathbb{E}\left[\left(\int_0^t\left(1-x_{i,s-}\right)\mathrm{d}\tilde{N}_{i,s}\right)^2+\left(\int_0^t x_{i,s-}\mathrm{d}\tilde{L}_{i,s}\right)^2\right]\\
&=\frac{1}{M^2}\sum_{i=1}^{M}\mathbb{E}\left[\left(\int_0^t\left(1-x_{i,s-}\right)^2 X_s^{(M)}g^+\left(X_s^{(M)}\right)\mathrm{d}s+\int_0^t\left(x_{i,s-}\right)^2\left(1-X_s^{(M)}\right)g^-\left(X_s^{(M)}\right)\mathrm{d}s\right)\right]\\
&=\frac{1}{M^2}\sum_{i=1}^{M}\mathbb{E}\left[\left(\int_0^t\left(1-x_{i,s-}\right)X_s^{(M)}g^+\left(X_s^{(M)}\right)\mathrm{d}s+\int_0^t x_{i,s-}\left(1-X_s^{(M)}\right)g^-\left(X_s^{(M)}\right)\mathrm{d}s\right)\right],\\
&=\frac{1}{M}\mathbb{E}\left[\left(\int_0^t\left(1-X_s^{(M)}\right)X_s^{(M)}g^+\left(X_s^{(M)}\right)\mathrm{d}s+\int_0^t X_s^{(M)}\left(1-X_s^{(M)}\right)g^-\left(X_s^{(M)}\right)\mathrm{d}s\right)\right]\\
&\leq\frac{C_g T}{M}
\end{aligned}
\tag{48}
$$

where $C_g=\max_{0\leq y\leq 1}\left\{\left(1-y\right)y\left(g^+\left(y\right)+g^-\left(y\right)\right)\right\}<+\infty$ is a nonnegative constant independent of $T$. Combining (47) and (48) yields

$$
\mathbb{E}\left[\left(\frac{1}{M}\sum_{i=1}^{M}\left\{\int_0^t\left(1-x_{i,s-}\right)\mathrm{d}\tilde{N}_{i,s}-\int_0^t x_{i,s-}\mathrm{d}\tilde{L}_{i,s}\right\}\right)^2\right]\leq\frac{C_g T}{M}\to 0 \text{ as } M\to+\infty,
\tag{49}
$$

proving (16). The second equality (17) is due to the first one combined with (15).

□

***Proof of Proposition 3***

The relationship (48) establishes that the process

$$
A_t^{(M)}\equiv\frac{1}{M}\sum_{i=1}^{M}\left\{\int_0^t\left(1-x_{i,s-}\right)\mathrm{d}\tilde{N}_{i,s}-\int_0^t x_{i,s-}\mathrm{d}\tilde{L}_{i,s}\right\},\quad t\in\left[0,T\right]
\tag{50}
$$

satisfies the Aldous criterion for the sufficiency of tightness (e.g., Proof of Theorem 3.1 in Bansaye and Méléard [66]); for any $\varepsilon>0$ and $\eta>0$, there exist some $M_0\in\mathbb{N}$ and $\delta>0$ such that

$$
\sup_{M\geq M_0}\sup_{S,S'}\mathbb{P}\left(\left|A_S^{(M)}-A_{S'}^{(M)}\right|>\varepsilon\right)\leq\eta,
\tag{51}
$$

where $S$ and $S'$ are any stopping times defined on the natural filtration generated by all $N_i$ and $L_i$ such that $S \le S' \le \min\{S+\delta, T\}$. Specifically, it follows that

$$\begin{aligned} \mathbb{E}\left[\left|A_S^{(M)} - A_{S'}^{(M)}\right|\right]^2 &\le \mathbb{E}\left[\left(A_S^{(M)} - A_{S'}^{(M)}\right)^2\right] \\ &= \mathbb{E}\left[\left(\frac{1}{M}\sum_{i=1}^{M}\left\{\int_S^{S'}\left(1-x_{i,s-}\right)\mathrm{d}\tilde{N}_{i,s} - \int_S^{S'} x_{i,s-}\mathrm{d}\tilde{L}_{i,s}\right\}\right)^2\right], \\ &\le \frac{C_g\delta}{M} \end{aligned} \tag{52}$$

where the last line tends to 0 uniformly for each $M \in \mathbb{N}$ as $\delta \to +0$, where the Markov inequality was used to derive the convergence in probability from that of the least-squares (e.g., Theorem 5.11 in Klenke [41]). These observations combined with **Proposition 2** demonstrate that the process

$$B_t^{(M)} \equiv \int_0^t \left(1 - X_s^{(M)}\right) X_s^{(M)} g\left(X_s^{(M)}\right) \mathrm{d}s \tag{53}$$

also satisfies the Aldous criterion. With the tightness of the processes $A^{(M)}$ and $B^{(M)}$, the tightness of the sequence of laws of $\left(X_t^{(M)}\right)_{0 \le t \le T}$ follows due to its strict boundedness.

Then, identifying the limit equation to be satisfied by the process $\hat{X}$ and following the methodology outlined on p. 22 in the Proof of Theorem 3.1 in Bansaye and Méléard [66] in conjunction with (15) and **Proposition 2**, the limit process $\hat{X}$ of $X^{(M)}$ satisfies (18), thereby completing the proof. Here, the tightness of $X^{(M)}$ combined with Prokhorov's theorem implies the relative compactness of the family of laws of $X^{(M)}$ ($M = 1,2,3,\ldots$) in the set of probability measures in the space $\mathbb{D}\left(\left[0,T\right],\left[0,1\right]\right)$, which leads to the limiting law.

□

***Proof of Proposition 4***

Fix $T > 0$ that will be chosen to be smaller if necessary. Let $C > 1$ denote a sufficiently large constant depending solely on $g^+, g^-$. First, the boundedness and contraction properties of $\mathbb{G}$ are established to demonstrate that equation (25) admits a unique solution in $\mathbb{L}_{1,T}$ for any $T > 0$. The continuity and continuous differentiability of this solution is then proven. Finally, the range of this solution is proven to be at most $\left[0,1\right]$, and hence, the solution also satisfies (21).

***Boundedness of*** $\mathbb{G}$

For any $\hat{x} \in \mathbb{L}_{1,T}$,

$$
\begin{aligned}
\left\|\mathbb{G}[\hat{x}]\right\|_T &= \left\| e^{-Rt}\hat{x}_0(R) + \int_0^t e^{-R(t-s)}\left\{ \tilde{X}_s g^+\left(\tilde{X}_s\right)\left(1-\hat{x}_s(R)\right) - \left(1-\tilde{X}_s\right) g^-\left(\tilde{X}_s\right)\tilde{x}_s(R) \right\} \mathrm{d}s \right\|_T \\
&\leq C\left( 1 + \left\| \int_0^t e^{-R(t-s)} \mathrm{d}s \right\|_T \right) \\
&= C(1+T)
\end{aligned}
\tag{54}
$$

and hence

$$
\left\|\mathbb{G}[\hat{x}]\right\|_T \leq C(1+T), \tag{55}
$$

showing that $\mathbb{G}$ maps $\mathbb{L}_{1,T}$ to $\mathbb{L}_{1,T}$ .

***Strict contraction of*** $\mathbb{G}$

For any $\hat{x},\hat{y} \in \mathbb{L}_{1,T}$ that share the common initial condition $\hat{x}_0$ (an analogous notation, i.e., (26), to $\hat{y}$ is used),

$$
\begin{aligned}
\left\|\mathbb{G}[\hat{x}] - \mathbb{G}[\hat{y}]\right\|_T &= \left\| \begin{matrix} \int_0^t e^{-R(t-s)}\left\{ \tilde{X}_s g^+\left(\tilde{X}_s\right)\left(1-\tilde{x}_s(R)\right) - \left(1-\tilde{X}_s\right) g^-\left(\tilde{X}_s\right)\tilde{x}_s(R) \right\} \mathrm{d}s \\ -\int_0^t e^{-R(t-s)}\left\{ \tilde{Y}_s g^+\left(\tilde{Y}_s\right)\left(1-\tilde{y}_s(R)\right) - \left(1-\tilde{Y}_s\right) g^-\left(\tilde{Y}_s\right)\tilde{y}_s(R) \right\} \mathrm{d}s \end{matrix} \right\|_T \\
&\leq \left\| \int_0^t e^{-R(t-s)}\left\{ \hat{X}_s g^+\left(\tilde{X}_s\right)\left(1-\tilde{x}_s(R)\right) \right\} \mathrm{d}s - \int_0^t e^{-R(t-s)}\left\{ \tilde{Y}_s g^+\left(\tilde{Y}_s\right)\left(1-\tilde{y}_s(R)\right) \right\} \mathrm{d}s \right\|_T \\
&\quad + \left\| \int_0^t e^{-R(t-s)}\left(1-\tilde{X}_s\right) g^-\left(\tilde{X}_s\right)\tilde{x}_s(R) \mathrm{d}s - \int_0^t e^{-R(t-s)}\left(1-\tilde{Y}_s\right) g^-\left(\tilde{Y}_s\right)\tilde{y}_s(R) \mathrm{d}s \right\|_T
\end{aligned}
\tag{56}
$$

The first term in the last line of (56) is evaluated as

$$
\begin{aligned}
&\left\| \int_0^t e^{-R(t-s)}\left\{ \tilde{X}_s g^+\left(\tilde{X}_s\right)\left(1-\tilde{x}_s(R)\right) \right\} \mathrm{d}s - \int_0^t e^{-R(t-s)}\left\{ \tilde{Y}_s g^+\left(\tilde{Y}_s\right)\left(1-\tilde{y}_s(R)\right) \right\} \mathrm{d}s \right\|_T \\
&= \left\| \begin{matrix} \int_0^t e^{-R(t-s)}\left\{ \tilde{X}_s g^+\left(\tilde{X}_s\right)\left(1-\tilde{x}_s(R)\right) \right\} \mathrm{d}s - \int_0^t e^{-R(t-s)}\left\{ \tilde{X}_s g^+\left(\tilde{X}_s\right)\left(1-\tilde{y}_s(R)\right) \right\} \mathrm{d}s \\ \int_0^t e^{-R(t-s)}\left\{ \tilde{X}_s g^+\left(\tilde{X}_s\right)\left(1-\tilde{y}_s(R)\right) \right\} \mathrm{d}s - \int_0^t e^{-R(t-s)}\left\{ \tilde{Y}_s g^+\left(\tilde{Y}_s\right)\left(1-\tilde{y}_s(R)\right) \right\} \mathrm{d}s \end{matrix} \right\|_T \\
&\leq C\left\| \int_0^t e^{-R(t-s)}\left|\tilde{y}_s(R) - \tilde{x}_s(R)\right| \mathrm{d}s \right\|_T + C\left\| \int_0^t e^{-R(t-s)}\left\{ \tilde{X}_s g^+\left(\tilde{X}_s\right) - \tilde{Y}_s g^+\left(\tilde{Y}_s\right) \right\} \mathrm{d}s \right\|_T \\
&\leq C\left\| \int_0^t e^{-R(t-s)}\left|\tilde{y}_s(R) - \tilde{x}_s(R)\right| \mathrm{d}s \right\|_T + C\left\| \int_0^t e^{-R(t-s)}\left|\tilde{X}_s - \tilde{Y}_s\right| \mathrm{d}s \right\|_T
\end{aligned}
\tag{57}
$$

It follows that

$$
\begin{aligned}
\left\| \int_0^t e^{-R(t-s)}\left|\tilde{y}_s(R) - \tilde{x}_s(R)\right| \mathrm{d}s \right\|_T &= \sup_{0\leq t\leq T}\left( \int_0^{+\infty}\int_0^t e^{-R(t-s)}\left|\tilde{y}_s(R) - \tilde{x}_s(R)\right| \mathrm{d}s F(\mathrm{d}R) \right) \\
&\leq \sup_{0\leq t\leq T}\left( \int_0^t\left( \int_0^{+\infty}\left|\tilde{y}_s(R) - \tilde{x}_s(R)\right| F(\mathrm{d}R) \right)\mathrm{d}s \right) \\
&\leq \sup_{0\leq t\leq T}\int_0^t \left\|\tilde{y} - \tilde{x}\right\|_T \mathrm{d}s \\
&= T\left\|\tilde{y} - \tilde{x}\right\|_T \\
&\leq T\left\|\hat{y} - \hat{x}\right\|_T
\end{aligned}
\tag{58}
$$

and

$$
\begin{aligned}
\left\| \int_0^t e^{-R(t-s)} \left| \tilde{X}_s - \tilde{Y}_s \right| \mathrm{d}s \right\|_T &\le \sup_{0 \le t \le T} \left( \int_0^t \left| \tilde{X}_s - \tilde{Y}_s \right| \mathrm{d}s \right) \\
&= \sup_{0 \le t \le T} \left( \int_0^t \left( \int_0^{+\infty} \left| \tilde{y}_s(R) - \tilde{x}_s(R) \right| F(\mathrm{d}R) \right) \mathrm{d}s \right). \\
&\le T \left\| \hat{y} - \hat{x} \right\|_T
\end{aligned}
\tag{59}
$$

Hence, by (57)–(59),

$$
\left\| \int_0^t e^{-R(t-s)} \left\{ \tilde{X}_s g^+\left(\tilde{X}_s\right)\left(1 - \tilde{x}_s(R)\right) \right\} \mathrm{d}s - \int_0^t e^{-R(t-s)} \left\{ \tilde{Y}_s g^+\left(\tilde{Y}_s\right)\left(1 - \tilde{y}_s(R)\right) \right\} \mathrm{d}s \right\|_T \le 2CT \left\| \hat{y} - \hat{x} \right\|_T . \tag{60}
$$

Considering the second term in the last line of (56), as in the discussion above, it follows that

$$
\begin{aligned}
&\left\| \int_0^t e^{-R(t-s)} \left(1 - \tilde{X}_s\right) g^-\left(\tilde{X}_s\right) \tilde{x}_s(R) \mathrm{d}s - \int_0^t e^{-R(t-s)} \left(1 - \tilde{Y}_s\right) g^-\left(\tilde{Y}_s\right) \tilde{y}_s(R) \mathrm{d}s \right\|_T \\
&= \left\| \begin{array}{l} \int_0^t e^{-R(t-s)} \left(1 - \tilde{X}_s\right) g^-\left(\tilde{X}_s\right) \left(\tilde{x}_s(R) - \tilde{y}_s(R)\right) \mathrm{d}s \\ + \int_0^t e^{-R(t-s)} \left\{ \left(1 - \tilde{X}_s\right) g^-\left(\tilde{X}_s\right) - \left(1 - \tilde{Y}_s\right) g^-\left(\tilde{Y}_s\right) \right\} \tilde{y}_s(R) \mathrm{d}s \end{array} \right\|_T \\
&\le C \left\| \int_0^t e^{-R(t-s)} \left| \tilde{x}_s(R) - \tilde{y}_s(R) \right| \mathrm{d}s \right\|_T + C \left\| \int_0^t e^{-R(t-s)} \left| \tilde{X}_s - \tilde{Y}_s \right| \mathrm{d}s \right\|_T
\end{aligned}
\tag{61}
$$

and hence

$$
\left\| \int_0^t e^{-R(t-s)} \left(1 - \tilde{X}_s\right) g^-\left(\tilde{X}_s\right) \tilde{x}_s(R) \mathrm{d}s - \int_0^t e^{-R(t-s)} \left(1 - \tilde{Y}_s\right) g^-\left(\tilde{Y}_s\right) \tilde{y}_s(R) \mathrm{d}s \right\|_T \le 2CT \left\| \hat{y} - \hat{x} \right\|_T . \tag{62}
$$

Consequently, it follows that

$$
\left\| \mathbb{G}[\hat{x}] - \mathbb{G}[\hat{y}] \right\|_T \le 4CT \left\| \hat{y} - \hat{x} \right\|_T . \tag{63}
$$

The mapping $\mathbb{G}$ is strictly contractive if $T < \dfrac{1}{4C}$. For such a small $T$, which can be selected irrespective of any $\hat{x}, \hat{y}$, the combination of boundedness and contraction results, along with the Banach fixed-point theorem, ensures that (25) admits a unique solution in $\mathbb{L}_{1,T}$. Since $T$ is independent of the solution itself, the above procedure can be repeated to the time intervals $[2T, 3T]$, $[3T, 4T]$,…, thereby continuing the solution globally in time.

***Continuity and smoothness of*** $\mathbb{G}[\hat{x}]$

The continuity of $\mathbb{G}[\hat{x}]$ is proven, where $\hat{x}$ denotes the unique solution to (25) derived above. Fix $T > 0$ and $R > 0$. Consider any $t, u \in [0, T]$ with $0 \le u \le t \le T$. With a constant $C' > 0$ not depending on $t, u$, it follows that

$$
\begin{aligned}
&\left|\mathbb{G}[\hat{x}](t,R)-\mathbb{G}[\hat{x}](u,R)\right| \\
&=\begin{vmatrix}
e^{-Rt}\hat{x}_0(R)+\int_0^t e^{-R(t-s)}\left\{\tilde{X}_s g^{+}\left(\tilde{X}_s\right)\left(1-\tilde{x}_s(R)\right)-\left(1-\tilde{X}_s\right)g^{-}\left(\tilde{X}_s\right)\tilde{x}_s(R)\right\}\mathrm{d}s \\
-e^{-Ru}\hat{x}_0(R)-\int_0^u e^{-R(u-s)}\left\{\tilde{X}_s g^{+}\left(\tilde{X}_s\right)\left(1-\tilde{x}_s(R)\right)-\left(1-\tilde{X}_s\right)g^{-}\left(\tilde{X}_s\right)\tilde{x}_s(R)\right\}\mathrm{d}s
\end{vmatrix} \\
&\leq\left|e^{-Rt}-e^{-Ru}\right|\hat{x}_0(R) \\
&\quad+\begin{vmatrix}
\int_0^t e^{-R(t-s)}\left\{\tilde{X}_s g^{+}\left(\tilde{X}_s\right)\left(1-\tilde{x}_s(R)\right)-\left(1-\tilde{X}_s\right)g^{-}\left(\tilde{X}_s\right)\tilde{x}_s(R)\right\}\mathrm{d}s \\
-\int_0^u e^{-R(t-s)}\left\{\tilde{X}_s g^{+}\left(\tilde{X}_s\right)\left(1-\tilde{x}_s(R)\right)-\left(1-\tilde{X}_s\right)g^{-}\left(\tilde{X}_s\right)\tilde{x}_s(R)\right\}\mathrm{d}s
\end{vmatrix} \\
&\quad+\begin{vmatrix}
\int_0^u e^{-R(t-s)}\left\{\tilde{X}_s g^{+}\left(\tilde{X}_s\right)\left(1-\tilde{x}_s(R)\right)-\left(1-\tilde{X}_s\right)g^{-}\left(\tilde{X}_s\right)\tilde{x}_s(R)\right\}\mathrm{d}s \\
-\int_0^u e^{-R(u-s)}\left\{\tilde{X}_s g^{+}\left(\tilde{X}_s\right)\left(1-\tilde{x}_s(R)\right)-\left(1-\tilde{X}_s\right)g^{-}\left(\tilde{X}_s\right)\tilde{x}_s(R)\right\}\mathrm{d}s
\end{vmatrix} \\
&\leq\left|e^{-Rt}-e^{-Ru}\right|\hat{x}_0(R)+C'\left(\int_u^t e^{-R(t-s)}\mathrm{d}s+\int_0^u\left|e^{-R(t-s)}-e^{-R(u-s)}\right|\mathrm{d}s\right)
\end{aligned}
\quad , \tag{64}
$$

confirming the continuity of $\mathbb{G}[\hat{x}]$ in time. It also follows that

$$
\left|\hat{x}_t(R)-\hat{x}_u(R)\right|\leq\left|e^{-Rt}-e^{-Ru}\right|\hat{x}_0(R)+C'\left(\int_u^t e^{-R(t-s)}\mathrm{d}s+\int_0^u\left|e^{-R(t-s)}-e^{-R(u-s)}\right|\mathrm{d}s\right) \tag{65}
$$

by the definition of $\hat{x}$. According to Proposition 2' in Zorich [85], $\mathbb{G}[\hat{x}]$ is a continuously differentiable function on $[0,T]$, owing to the boundedness and continuity of $\hat{x}$. Since $T>0$ is arbitrary, it follows that the solution $\hat{x}$ to (25) is continuously differentiable for any $t>0$; therefore, this solution satisfies (27) pointwise.

***Boundedness of*** $\hat{x}$

The equation (27) is rewritten as

$$
\frac{\mathrm{d}\hat{x}_t(R)}{\mathrm{d}t}=\varpi\left(R,\hat{x}_t(R)\right),\quad 0<t\leq T\ ,\quad R>0\ . \tag{66}
$$

The initial condition $\hat{x}_0$ is bounded between 0 and 1. Assume that $\hat{x}_t(R)=0$ for some $(t,R)$. Then, it follows that

$$
\varpi\left(R,\hat{x}_t(R)\right)=\tilde{X}_t g^{+}\left(\tilde{X}_t\right)\geq 0\ . \tag{67}
$$

Similarly, assume that $\hat{x}_t(R)=1$ at some $(t,R)$. Then, it follows that

$$
\varpi\left(R,\hat{x}_t(R)\right)=-\left\{R+\left(1-\tilde{X}_t\right)g^{-}\left(\tilde{X}_t\right)\right\}\leq 0\ . \tag{68}
$$

Consequently, $\hat{x}$ must remain bounded between 0 and 1 and satisfy (21) pointwise.

□

***Proof of Proposition 5***

By **Proposition 4**, for any $t,R>0$, the regularized model (25) yields

$$
\hat{x}_t(R)=\mathbb{G}[\hat{x}](t,R)\quad\text{and}\quad\hat{y}_t(R)=\mathbb{G}[\hat{y}](t,R)\ . \tag{69}
$$

By the strict boundedness of the solutions $\hat{x}, \hat{y}$ (i.e., which are bounded between 0 and 1) and the Lipschitz continuity of the coefficients in the integrands of $\mathbb{G}$, it follows that

$$
\begin{aligned}
&\left|\hat{x}_t(R)-\hat{y}_t(R)\right| \\
&=\left|\mathbb{G}[\hat{x}](t,R)-\mathbb{G}[\hat{y}](t,R)\right| \\
&\leq\left|\hat{x}_0(R)-\hat{y}_0(R)\right|e^{-Rt}+\left|\begin{array}{l}\int_0^t e^{-R(t-s)}\left\{\tilde{X}_s g^{+}\left(\tilde{X}_s\right)\left(1-\tilde{x}_s(R)\right)-\left(1-\tilde{X}_s\right)g^{-}\left(\tilde{X}_s\right)\tilde{x}_s(R)\right\}\mathrm{d}s \\ -\int_0^t e^{-R(t-s)}\left\{\tilde{X}_s g^{+}\left(\tilde{X}_s\right)\left(1-\tilde{y}_s(R)\right)-\left(1-\tilde{X}_s\right)g^{+}\left(\tilde{X}_s\right)\tilde{y}_s(R)\right\}\mathrm{d}s\end{array}\right| \\
&+\left|\begin{array}{l}+\int_0^t e^{-R(t-s)}\left\{\tilde{X}_s g^{+}\left(\tilde{X}_s\right)\left(1-\tilde{y}_s(R)\right)-\left(1-\tilde{X}_s\right)g^{+}\left(\tilde{X}_s\right)\tilde{y}_s(R)\right\}\mathrm{d}s \\ -\int_0^t e^{-R(t-s)}\left\{\tilde{Y}_s g^{+}\left(\tilde{Y}_s\right)\left(1-\tilde{y}_s(R)\right)-\left(1-\tilde{Y}_s\right)g^{-}\left(\tilde{Y}_s\right)\tilde{y}_s(R)\right\}\mathrm{d}s\end{array}\right| \\
&\leq\left|\hat{x}_0(R)-\hat{y}_0(R)\right|e^{-Rt}+\int_0^t e^{-R(t-s)}\left(\tilde{X}_s g^{+}\left(\tilde{X}_s\right)+\left(1-\tilde{X}_s\right)g^{+}\left(\tilde{X}_s\right)\right)\left|\tilde{x}_s(R)-\tilde{y}_s(R)\right|\mathrm{d}s \\
&+\int_0^t e^{-R(t-s)}\left(1-\tilde{y}_s(R)\right)\left|\tilde{X}_s g^{+}\left(\tilde{X}_s\right)-\tilde{Y}_s g^{+}\left(\tilde{Y}_s\right)\right|\mathrm{d}s \\
&+\int_0^t e^{-R(t-s)}\tilde{y}_s(R)\left|\left(1-\tilde{X}_s\right)g^{+}\left(\tilde{X}_s\right)-\left(1-\tilde{Y}_s\right)g^{-}\left(\tilde{Y}_s\right)\right|\mathrm{d}s \\
&\leq\left|\hat{x}_0(R)-\hat{y}_0(R)\right|+C\int_0^t\left|\tilde{X}_s-\tilde{Y}_s\right|\mathrm{d}s
\end{aligned}
\tag{70}
$$

where (similar notation applies to $x$-based ones)

$$
\tilde{y}_s(R)=\max\left\{0,\min\left\{1,\hat{y}_s(R)\right\}\right\} \quad \text{and} \quad \tilde{Y}_s=\int_0^{+\infty}\tilde{y}_s(R)F_y(\mathrm{d}R). \tag{71}
$$

It follows that

$$
\begin{aligned}
\int_0^t\left|\tilde{X}_s-\tilde{Y}_s\right|\mathrm{d}s &=\int_0^t\left|\int_0^{+\infty}\tilde{x}_s(R)F_x(\mathrm{d}R)-\int_0^{+\infty}\tilde{y}_s(R)F_y(\mathrm{d}R)\right|\mathrm{d}s \\
&\leq\int_0^t\left|\int_0^{+\infty}\tilde{x}_s(R)\left(F_x(\mathrm{d}R)-F_y(\mathrm{d}R)\right)\right|\mathrm{d}s+\int_0^t\int_0^{+\infty}\left|\tilde{x}_s(R)-\tilde{y}_s(R)\right|F_y(\mathrm{d}R)\mathrm{d}s . \\
&\leq 2\int_0^t\left\|F_x-F_y\right\|_{\mathrm{TV}}\mathrm{d}s+\int_0^t\int_0^{+\infty}\left|\tilde{x}_s(R)-\tilde{y}_s(R)\right|F_y(\mathrm{d}R)\mathrm{d}s
\end{aligned}
\tag{72}
$$

By (70)–(72), with a sufficiently large $C \geq 2$,

$$
\left|\hat{x}_t(R)-\hat{y}_t(R)\right|\leq\left|\hat{x}_0(R)-\hat{y}_0(R)\right|+C\left(\int_0^t\left\|F_x-F_y\right\|_{\mathrm{TV}}\mathrm{d}s+\int_0^t\int_0^{+\infty}\left|\tilde{x}(Q)-\tilde{y}(Q)\right|F_y(\mathrm{d}Q)\mathrm{d}s\right). \tag{73}
$$

The right-hand side of (73) is further evaluated as

$$
\begin{aligned}
&\left|\hat{x}_0(R)-\hat{y}_0(R)\right|+C\left(\int_0^t\left\|F_x-F_y\right\|_{\mathrm{TV}}\mathrm{d}s+\int_0^t\int_0^{+\infty}\left|\tilde{x}_s(Q)-\tilde{y}_s(Q)\right|F_y(\mathrm{d}Q)\mathrm{d}s\right) \\
&\leq\sup_{R>0}\left|\hat{x}_0(R)-\hat{y}_0(R)\right|+C\int_0^t\left(\left\|F_x-F_y\right\|_{\mathrm{TV}}+\sup_{R>0}\left|\hat{x}_s(R)-\hat{y}_s(R)\right|\right)\mathrm{d}s
\end{aligned}
, \tag{74}
$$

yielding

$$
\left|\hat{x}_t(R)-\hat{y}_t(R)\right|\leq\sup_{R>0}\left|\hat{x}_0(R)-\hat{y}_0(R)\right|+C\int_0^t\left(\left\|F_x-F_y\right\|_{\mathrm{TV}}+\sup_{R>0}\left|\hat{x}_s(R)-\hat{y}_s(R)\right|\right)\mathrm{d}s, \tag{75}
$$

and hence

$$
\sup_{R>0}\left|\hat{x}_t(R)-\hat{y}_t(R)\right|\leq\sup_{R>0}\left|\hat{x}_0(R)-\hat{y}_0(R)\right|+C\int_0^t\left(\left\|F_x-F_y\right\|_{\mathrm{TV}}+\sup_{R>0}\left|\hat{x}_s(R)-\hat{y}_s(R)\right|\right)\mathrm{d}s . \tag{76}
$$

Applying Gronwall's inequality (e.g., Gronwall [86]) to (76) and again selecting a suitably large $C > 0$ if necessary yields (28).

□

### A.3 A Lemma

***Lemma 1***

*The stochastic system (12) admits at most one càdlàg pathwise solution that is bounded between 0 and 1.*

***Proof of Lemma 1***

The system reads, for each $i = 1,2,3,...,M$,

$$\begin{aligned} x_{i,t} = x_{i,0} &+ \int_0^t \left(1 - x_{i,s}\right) X_s^{(M)} g^+\left(X_s^{(M)}\right) \mathrm{d}s - \int_0^t x_{i,s} \left(1 - X_s^{(M)}\right) g^-\left(X_s^{(M)}\right) \mathrm{d}s \\ &+ \int_0^t \left(1 - x_{i,s-}\right) \mathrm{d}\tilde{N}_{i,s} - \int_0^t x_{i,s-} \mathrm{d}\tilde{L}_{i,s} \end{aligned} . \tag{77}$$

Throughout the proof, $C > 0$ denotes a constant whose value may vary line by line.

Fix $t > 0$ and $i \in \{1,2,3,...,M\}$. Assume that there exist two pathwise solutions $y_i$ and $z_i$ (bounded almost surely between 0 and 1) to (77) with a common initial condition $x_{i,0}$. It follows that (by using the representations (13))

$$\begin{aligned} &\left|y_{i,t} - z_{i,t}\right|^2 \\ &= \left\{ \begin{aligned} &\int_0^t \left(1 - y_{i,s}\right) Y_s^{(M)} g^+\left(Y_s^{(M)}\right) \mathrm{d}s - \int_0^t \left(1 - z_{i,s}\right) Z_s^{(M)} g^+\left(Z_s^{(M)}\right) \mathrm{d}s \\ &-\left( \int_0^t y_{i,s} \left(1 - Y_s^{(M)}\right) g^-\left(Y_s^{(M)}\right) \mathrm{d}s - \int_0^t z_{i,s} \left(1 - Z_s^{(M)}\right) g^-\left(Z_s^{(M)}\right) \mathrm{d}s \right) \\ &+\int_0^t \left(1 - y_{i,s-}\right) \int_{u=0}^{u=Y_{s-}^{(M)} g^+\left(Y_{s-}^{(M)}\right)} \int_{z=0}^{z=1} \tilde{\bar{N}}_i\left(\mathrm{d}u\mathrm{d}z\mathrm{d}s\right) - \int_0^t \left(1 - z_{i,s-}\right) \int_{u=0}^{u=Z_{s-}^{(M)} g^+\left(Z_{s-}^{(M)}\right)} \int_{z=0}^{z=1} \tilde{\bar{N}}_i\left(\mathrm{d}u\mathrm{d}z\mathrm{d}s\right) \\ &-\left( \int_0^t y_{i,s-} \int_{u=0}^{u=\left(1-Y_{s-}^{(M)}\right) g^-\left(Y_{s-}^{(M)}\right)} \int_{z=0}^{z=1} \tilde{\bar{L}}_i\left(\mathrm{d}u\mathrm{d}z\mathrm{d}s\right) - \int_0^t z_{i,s-} \int_{u=0}^{u=\left(1-Z_{s-}^{(M)}\right) g^-\left(Z_{s-}^{(M)}\right)} \int_{z=0}^{z=1} \tilde{\bar{L}}_i\left(\mathrm{d}u\mathrm{d}z\mathrm{d}s\right) \right) \end{aligned} \right\}^2 \\ &\leq 4\left\{ \int_0^t \left(1 - y_{i,s-}\right) Y_s^{(M)} g^+\left(Y_s^{(M)}\right) \mathrm{d}s - \int_0^t \left(1 - z_{i,s-}\right) Z_s^{(M)} g^+\left(Z_s^{(M)}\right) \mathrm{d}s \right\}^2 \\ &+4\left\{ \int_0^t \left(1 - y_{i,s-}\right) Y_s^{(M)} g^+\left(Y_s^{(M)}\right) \mathrm{d}s - \int_0^t \left(1 - z_{i,s-}\right) Z_s^{(M)} g^+\left(Z_s^{(M)}\right) \mathrm{d}s \right\}^2 \\ &+4\left\{ \int_0^t \left(1 - y_{i,s-}\right) \int_{u=0}^{u=Y_{s-}^{(M)} g^+\left(Y_{s-}^{(M)}\right)} \int_{z=0}^{z=1} \tilde{\bar{N}}_i\left(\mathrm{d}u\mathrm{d}z\mathrm{d}s\right) - \int_0^t \left(1 - z_{i,s-}\right) \int_{u=0}^{u=Z_{s-}^{(M)} g^+\left(Z_{s-}^{(M)}\right)} \int_{z=0}^{z=1} \tilde{\bar{N}}_i\left(\mathrm{d}u\mathrm{d}z\mathrm{d}s\right) \right\}^2 \\ &+4\left\{ \int_0^t y_{i,s-} \int_{u=0}^{u=\left(1-Y_{s-}^{(M)}\right) g^-\left(Y_{s-}^{(M)}\right)} \int_{z=0}^{z=1} \tilde{\bar{L}}_i\left(\mathrm{d}u\mathrm{d}z\mathrm{d}s\right) - \int_0^t z_{i,s-} \int_{u=0}^{u=\left(1-Z_{s-}^{(M)}\right) g^-\left(Z_{s-}^{(M)}\right)} \int_{z=0}^{z=1} \tilde{\bar{L}}_i\left(\mathrm{d}u\mathrm{d}z\mathrm{d}s\right) \right\}^2 \end{aligned} . \tag{78}$$

Here, $\tilde{\bar{N}}_i$ and $\tilde{\bar{L}}_i$ are compensated versions of $\bar{N}_i$ and $\bar{L}_i$, respectively.

Each term in (78) is evaluated as follows. The first term in (78) is evaluated as

$$
\begin{aligned}
&\left\{\int_0^t\left(1-y_{i,s}\right)Y_s^{(M)}g^+\left(Y_s^{(M)}\right)\mathrm{d}s-\int_0^t\left(1-z_{i,s}\right)Z_s^{(M)}g^+\left(Z_s^{(M)}\right)\mathrm{d}s\right\}^2\\
&=\left\{\int_0^t\left(1-y_{i,s}\right)\left(Y_s^{(M)}g^+\left(Y_s^{(M)}\right)-Z_s^{(M)}g^+\left(Z_s^{(M)}\right)\right)\mathrm{d}s+\int_0^t\left(z_{i,s}-y_{i,s}\right)Z_s^{(M)}g^+\left(Z_s^{(M)}\right)\mathrm{d}s\right\}^2\\
&\leq 2\left\{\left(\int_0^t\left(1-y_{i,s}\right)\left(Y_s^{(M)}g^+\left(Y_s^{(M)}\right)-Z_s^{(M)}g^+\left(Z_s^{(M)}\right)\right)\mathrm{d}s\right)^2+\left(\int_0^t\left|z_{i,s}-y_{i,s}\right|Z_s^{(M)}g^+\left(Z_s^{(M)}\right)\mathrm{d}s\right)^2\right\}\\
&\leq C\left\{\left(\int_0^t\left|Y_s^{(M)}g^+\left(Y_s^{(M)}\right)-Z_s^{(M)}g^+\left(Z_s^{(M)}\right)\right|\mathrm{d}s\right)^2+\left(\int_0^t\left|z_{i,s}-y_{i,s}\right|\mathrm{d}s\right)^2\right\}\\
&\leq C\left\{\left(\int_0^t\left|Y_s^{(M)}-Z_s^{(M)}\right|\mathrm{d}s\right)^2+\left(\int_0^t\left|z_{i,s}-y_{i,s}\right|\mathrm{d}s\right)^2\right\}\\
&\leq CT\left\{\int_0^t\left|Y_s^{(M)}-Z_s^{(M)}\right|^2\mathrm{d}s+\int_0^t\left|z_{i,s}-y_{i,s}\right|^2\mathrm{d}s\right\}
\end{aligned}
\tag{79}
$$

Similarly, the second term in (78) is evaluated as

$$
\begin{aligned}
&\left\{\int_0^t\left(1-y_{i,s}\right)Y_s^{(M)}g^+\left(Y_s^{(M)}\right)\mathrm{d}s-\int_0^t\left(1-z_{i,s}\right)Z_s^{(M)}g^+\left(Z_s^{(M)}\right)\mathrm{d}s\right\}^2\\
&\leq CT\left\{\int_0^t\left|Y_s^{(M)}-Z_s^{(M)}\right|^2\mathrm{d}s+\int_0^t\left|z_{i,s}-y_{i,s}\right|^2\mathrm{d}s\right\}
\end{aligned}
\tag{80}
$$

The third term in (78) is evaluated as (one can apply the same method to the fourth term)

$$
\begin{aligned}
&\left\{\int_0^t\left(1-y_{i,s-}\right)\int_{u=0}^{u=Y_{s-}^{(M)}g^+\left(Y_{s-}^{(M)}\right)}\int_{z=0}^{z=1}\tilde{\tilde{N}}_i\left(\mathrm{d}u\mathrm{d}z\mathrm{d}s\right)-\int_0^t\left(1-z_{i,s-}\right)\int_{u=0}^{u=Z_{s-}^{(M)}g^+\left(Z_{s-}^{(M)}\right)}\int_{z=0}^{z=1}\tilde{\tilde{N}}_i\left(\mathrm{d}u\mathrm{d}z\mathrm{d}s\right)\right\}^2\\
&=\left\{\begin{array}{l}\displaystyle\int_0^t\left(1-y_{i,s-}\right)\int_{u=0}^{u=Y_{s-}^{(M)}g^+\left(Y_{s-}^{(M)}\right)}\int_{z=0}^{z=1}\tilde{\tilde{N}}_i\left(\mathrm{d}u\mathrm{d}z\mathrm{d}s\right)-\int_0^t\left(1-y_{i,s-}\right)\int_{u=0}^{u=Z_{s-}^{(M)}g^+\left(Z_{s-}^{(M)}\right)}\int_{z=0}^{z=1}\tilde{\tilde{N}}_i\left(\mathrm{d}u\mathrm{d}z\mathrm{d}s\right)\\ \displaystyle+\int_0^t\left(1-y_{i,s-}\right)\int_{u=0}^{u=Z_{s-}^{(M)}g^+\left(Z_{s-}^{(M)}\right)}\int_{z=0}^{z=1}\tilde{\tilde{N}}_i\left(\mathrm{d}u\mathrm{d}z\mathrm{d}s\right)-\int_0^t\left(1-z_{i,s-}\right)\int_{u=0}^{u=Z_{s-}^{(M)}g^+\left(Z_{s-}^{(M)}\right)}\int_{z=0}^{z=1}\tilde{\tilde{N}}_i\left(\mathrm{d}u\mathrm{d}z\mathrm{d}s\right)\end{array}\right\}^2\\
&\leq 2\left\{\int_0^t\left(1-y_{i,s-}\right)\int_{u=Z_{s-}^{(M)}g^+\left(Z_{s-}^{(M)}\right)}^{u=Y_{s-}^{(M)}g^+\left(Y_{s-}^{(M)}\right)}\int_{z=0}^{z=1}\tilde{\tilde{N}}_i\left(\mathrm{d}u\mathrm{d}z\mathrm{d}s\right)\right\}^2+2\left\{\int_0^t\left|z_{i,s-}-y_{i,s-}\right|\int_{u=0}^{u=Z_{s-}^{(M)}g^+\left(Z_{s-}^{(M)}\right)}\int_{z=0}^{z=1}\tilde{\tilde{N}}_i\left(\mathrm{d}u\mathrm{d}z\mathrm{d}s\right)\right\}^2\\
&\leq 2\left\{\int_0^t\int_{u=Z_{s-}^{(M)}g^+\left(Z_{s-}^{(M)}\right)}^{u=Y_{s-}^{(M)}g^+\left(Y_{s-}^{(M)}\right)}\int_{z=0}^{z=1}\tilde{\tilde{N}}_i\left(\mathrm{d}u\mathrm{d}z\mathrm{d}s\right)\right\}^2+2\left\{\int_0^t\left|z_{i,s-}-y_{i,s-}\right|\int_{u=0}^{u=Z_{s-}^{(M)}g^+\left(Z_{s-}^{(M)}\right)}\int_{z=0}^{z=1}\tilde{\tilde{N}}_i\left(\mathrm{d}u\mathrm{d}z\mathrm{d}s\right)\right\}^2
\end{aligned}
\tag{81}
$$

Each term in the last line of (81) is estimated as shown below.

The first term in the last line of (82), which is a square of a martingale, is evaluated as follows:

$$
\begin{aligned}
\mathbb{E}\left[\sup_{t\leq T}\left\{\int_0^t\int_{u=Z_{s-}^{(M)}g^+\left(Z_{s-}^{(M)}\right)}^{u=Y_{s-}^{(M)}g^+\left(Y_{s-}^{(M)}\right)}\int_{z=0}^{z=1}\tilde{\tilde{N}}_i\left(\mathrm{d}u\mathrm{d}z\mathrm{d}s\right)\right\}^2\right]&\leq C\mathbb{E}\left[\int_0^t\left|Y_s^{(M)}g^+\left(Y_s^{(M)}\right)-Z_s^{(M)}g^+\left(Z_s^{(M)}\right)\right|^2\mathrm{d}s\right]\\
&\leq C\mathbb{E}\left[\int_0^t\left|Y_s^{(M)}-Z_s^{(M)}\right|^2\mathrm{d}s\right]
\end{aligned}
\tag{82}
$$

Doob's inequality was used in the first line of (82) (e.g., 3.8 Theorem (iv) in Karatzas and Shreve [87] for $p=2$). The second term in the last line of (82) is evaluated as

$$
\begin{aligned}
\mathbb{E}\left[\left\{\int_0^t\int_{u=0}^{u=Z_{s-}^{(M)}g^+\left(Z_{s-}^{(M)}\right)}\int_{z=0}^{z=1}\left|z_{i,s-}-y_{i,s-}\right|\tilde{\tilde{N}}_i\left(\mathrm{d}u\mathrm{d}z\mathrm{d}s\right)\right\}^2\right]&\leq C\mathbb{E}\left[\int_0^t\left|z_{i,s}-y_{i,s}\right|^2Z_s^{(M)}g^+\left(Z_s^{(M)}\right)\mathrm{d}s\right]\\
&\leq C\mathbb{E}\left[\int_0^t\left|z_{i,s}-y_{i,s}\right|^2\mathrm{d}s\right]
\end{aligned}
\tag{83}
$$

Consequently, it follows that

$$
\begin{aligned}
&\mathbb{E}\left[\sup_{t\leq T}\left\{\int_0^t\left(1-y_{i,s-}\right)\int_{u=0}^{u=Y_{t-}^{(M)}g^+\left(Y_-^{(M)}\right)}\int_{z=0}^{z=1}\tilde{\tilde{N}}_i\left(\mathrm{d}u\mathrm{d}z\mathrm{d}s\right)-\int_0^t\left(1-z_{i,s-}\right)\int_{u=0}^{u=Z_{t-}^{(M)}g^+\left(Z_-^{(M)}\right)}\int_{z=0}^{z=1}\tilde{\tilde{N}}_i\left(\mathrm{d}u\mathrm{d}z\mathrm{d}s\right)\right\}^2\right] \\
&\leq C\left(\mathbb{E}\left[\int_0^t\left|Y_s^{(M)}-Z_s^{(M)}\right|^2\mathrm{d}s\right]+\mathbb{E}\left[\int_0^t\left|z_{i,s}-y_{i,s}\right|^2\mathrm{d}s\right]\right)
\end{aligned}. \tag{84}
$$

Set $U_{i,t}=\sup_{s\leq t}\left|y_{i,s}-z_{i,s}\right|^2$. Summarizing (78)-(84) yields

$$
\begin{aligned}
\mathbb{E}\left[U_{i,t}\right]&\leq C\left(T+1\right)\left(\mathbb{E}\left[\int_0^t\left|Y_s^{(M)}-Z_s^{(M)}\right|^2\mathrm{d}s\right]+\mathbb{E}\left[\int_0^t\left|z_{i,s}-y_{i,s}\right|^2\mathrm{d}s\right]\right) \\
&\leq C\left(T+1\right)\left(\mathbb{E}\left[\int_0^t\sup_{r\leq s}\left|Y_r^{(M)}-Z_r^{(M)}\right|^2\mathrm{d}s\right]+\mathbb{E}\left[\int_0^t U_{i,s}\mathrm{d}s\right]\right)
\end{aligned}. \tag{85}
$$

The integrand of the first integral in the last line of (85) is evaluated as follows:

$$
\sup_{r\leq s}\left|Y_r^{(M)}-Z_r^{(M)}\right|^2\leq\frac{C}{M}\sum_{j=1}^{M}U_{j,s}\ , \tag{86}
$$

which combined with (85) yields

$$
\mathbb{E}\left[U_{i,t}\right]\leq C\left(T+1\right)\left(\frac{1}{M}\sum_{j=1}^{M}\int_0^t\mathbb{E}\left[U_{j,s}\right]\mathrm{d}s+\int_0^t\mathbb{E}\left[U_{i,s}\right]\mathrm{d}s\right). \tag{87}
$$

Here, the exchange of the orders of expectation and integration is possible because each $U_{j,s}$ is bounded. Now, the estimate (87) for all $i=1,2,3,...,M$ follows, which combined with the Gronwall's inequality (e.g., Gronwall [86]) yields $\mathbb{E}\left[U_{i,t}\right]=0$ ($i=1,2,3,...,M$), and hence the uniqueness holds true.

□